%% file: decentr-revised-arXiv.tex
\documentclass[10pt,twocolumn,twoside]{IEEEtran}

\usepackage{etex}

\IEEEoverridecommandlockouts                              

\usepackage{psfrag}

\usepackage{latexsym}
\usepackage{graphicx}
\usepackage{float}

\usepackage{colortbl}
\usepackage{fancyhdr}
\usepackage{subfig}
\usepackage{hyperref}
\usepackage{enumerate}
\usepackage{multicol}
\usepackage{float}
\usepackage{cite}
\usepackage{mathtools, cuted}

\usepackage{graphicx,color}
\usepackage{amsmath,amssymb,amsfonts,booktabs, mathrsfs}
\usepackage{bbm,dsfont}
\usepackage[usenames,dvipsnames,svgnames,table]{xcolor}
\usepackage{tikz}
\usepackage{mhchem}
\usetikzlibrary{calc,arrows,automata,positioning,shapes,patterns}
\usepackage{pgfplots}
\usepgfplotslibrary{patchplots}

\newcommand{\tc}{\textcolor}

\newcommand{\card} {\mathrm{card}\,}
\newcommand{\sparse} {\mathrm{sp}\,}

\DeclareMathOperator{\blkdiag}{blkdiag}

\DeclareMathOperator*{\argmin}{argmin}
\DeclareMathOperator*{\minimize}{minimize}

\DeclareMathOperator*{\subject}{subject~to}

\newcommand{\mre}{\mathrm{e}}
\newcommand{\mrd}{\mathrm{d}}
\newcommand{\bbR}{\mathbb{R}}
\newcommand{\DefinedAs}[0]{\mathrel{\mathop:}=}

\definecolor{matlabblue}{rgb}{0   0.4470   0.7410}
\definecolor{matlabred}{rgb}{0.8500    0.3250    0.0980}
\definecolor{matlaborange}{rgb}{0.9290    0.6940    0.1250}
\newcommand{\vsp}{\vspace*{0.0cm}}
\newcommand{\vspa}{\vspace*{0.1cm}}

\newcommand{\cS}{{\cal S}}
\newcommand{\cC}{{\cal C}}
\newcommand{\Acl}{A_{\mathrm{cl}}}
\newcommand{\Jlb}{J^{\mathrm{lb}}}
\newcommand{\prox}{\mathbf{prox}}
\newcommand{\cG}{{\cal G}}
\newcommand{\one}{\mathds{1}}
\input{commands}

\usepackage{soul}
\usepackage[ruled]{algorithm2e}

\newcommand{\M}{\mathbb{M}}

\definecolor{cytoyellow}{RGB}{255, 204, 0}
\definecolor{cytoblue}{RGB}{0, 0, 204}
\definecolor{cytoblue2}{RGB}{51, 51, 255}
\definecolor{cytored}{RGB}{152, 0, 0}

\tikzset{
  treenode/.style = {align=center, inner sep=0pt, text centered,
    font=\sffamily},
  arn_n/.style = {treenode, circle, white, font=\sffamily\bfseries, draw=black,
    fill=black, text width=1.5em},
  arn_r/.style = {treenode, circle, white, font=\sffamily\bfseries, draw=red,
    fill=red, text width=1.5em},
  arn_x/.style = {treenode, rectangle, draw=black,
    minimum width=0.5em, minimum height=0.5em}
}

\begin{document}

\title{\Large \bf Structured decentralized control of positive systems with applications to
\\[0.1cm] combination drug therapy and leader selection in directed networks}
\author{Neil K.\ Dhingra, Marcello Colombino, and Mihailo R.\ Jovanovi\'c
\thanks{Financial support from the National Science Foundation under Awards ECCS-1739210 and CNS-1544887 and the Air Force Office of Scientific Research under Award FA9550-16-1-0009 is gratefully acknowledged.}
\thanks{N.\ K.\ Dhingra is with the Department of Electrical and Computer Engineering, University of Minnesota, M.\ Colombino is with the Automatic Control Laboratory, ETH Z\"urich, and M.\ R.\ Jovanovi\'c is with the Ming Hsieh Department of Electrical Engineering, University of Southern California. E-mails: dhin0008@umn.edu, mcolombi@control.ee.ethz.ch, mihailo@usc.edu.}}

\maketitle
\thispagestyle{empty}
\pagestyle{empty}

\begin{abstract}
We study a class of structured optimal control problems in which the main diagonal of the dynamic matrix is a linear function of the design variable. While such problems are in general challenging and nonconvex, for positive systems we prove convexity of the $\cH_2$ and $\cH_\infty$ optimal control formulations which allow for arbitrary convex constraints and regularization of the control input. Moreover, we establish differentiability of the $\Hinf$ norm when the graph associated with the dynamical generator is weakly connected and develop a customized algorithm for computing the optimal solution even in the absence of differentiability. We apply our results to the problems of leader selection in directed consensus networks and combination drug therapy for HIV treatment. In the context of leader selection, we address the combinatorial challenge by deriving upper and lower bounds on optimal performance. For combination drug therapy, we develop a customized subgradient method for efficient treatment of diseases whose mutation patterns are not connected.
\end{abstract}
	
	\vspace*{-3ex}
\section{Introduction}

Modern applications require structured controllers that cannot be designed using traditional approaches. Except in special cases, e.g., in funnel causal and quadratically invariant systems~\cite{bamvou05,rotlalTAC06}, and in system level synthesis approach~\cite{wang2016system} in which spatial and temporal sparsity constraints are imposed on the closed-loop response, posing optimal control problems in coordinates that preserve {\em controller structure\/} compromises convexity of the performance metric.  In this paper, we study the structured decentralized control of positive systems. While structured decentralized control is challenging in general, we show that, for positive systems, the convexity of the $\cH_2$ and $\cH_\infty$ optimal control formulations is not lost by imposing structural constraints. We also derive a graph theoretic condition that guarantees continuous differentiability of the $\cH_\infty$ performance metric and develop techniques to address combination drug therapy design and leader selection in directed \mbox{consensus networks.}

Positive systems arise in the modeling of systems with inherently nonnegative variables (e.g., probabilities, concentrations, and densities). Such systems have nonnegative states and outputs for nonnegative initial conditions and inputs~\cite{farrin11}. In this work, we examine models of HIV mutation~\cite{hernandez2011discrete,herestmid13,jonmatmur13,jonranmur13,jonmatmur14,jon16} and consensus networks~\cite{patbam10,linfarjovTAC14leaderselection,fitch2013information,fitch2014joint,clabuspoo14,clark2016submodularity}, where positivity comes from nonnegativity of populations and the structure of the underlying dynamics, respectively. Decentralized control, in which only the diagonal of the dynamical generator may be modified, is a suitable paradigm for modeling the effect of drugs on HIV~\cite{hernandez2011discrete} and the influence of leaders on the dynamics of leader-follower consensus networks~\cite{linfarjovTAC14leaderselection}. In these applications, the structure of decentralized control is important for capturing the efficacy of the drugs on different HIV mutants and influence of noise on the quality of consensus, respectively. This model can also be used to study chemical reaction networks and transportation networks.
	
The mathematical properties of positive systems can be exploited for efficient or structured controller design. In~\cite{tanlan11}, the authors show that the KYP lemma greatly simplifies for positive systems, thereby enabling decentralized $\mathcal H_\infty$ synthesis via Semidefinite Programming (SDP). In~\cite{ran15}, it is shown that a static output-feedback problem can be solved via Linear Programming (LP) for a class of positive systems. In~\cite{bri13,ebipeaarz11}, the authors develop necessary and sufficient conditions for robust stability of positive systems with respect to induced $\cL_1$--$\cL_\infty$ norm-bounded perturbations. In~\cite{ColSmi:2014:IFA_4769,ColSmi:2016:IFA_5242}, it is shown that the structured singular value is equal to its convex upper bound for positive systems so assessing robust stability with respect to induced $\cL_2$ norm-bounded perturbations \mbox{becomes tractable.}

It has been recently shown that the design of unconstrained decentralized controllers for positive systems can be cast as a convex problem~\cite{ran15,ran16}. However, since structural constraints cannot be handled by the LP or LMI approaches of~\cite{ran15,ran16}, references~\cite{jonranmur13,jonmatmur13,jonmatmur14,jon16} design $\cL_1$ and $\cH_\infty$ controllers that satisfy such constraints but achieve suboptimal performance. Furthermore, convexity of the weighted $\cL_1$ norm for structured decentralized control of positive systems was established in~\cite{colaneri2014convexity,ranber14} and optimized switching strategies were considered in~\cite{herestcol14,blacolval15,colmidbla16}. 

The paper is organized as follows. In Section~\ref{sec:problem:statement}, we formulate the regularized optimal control problem for a class positive systems. In Section~\ref{sec.conv}, we establish convexity of both $\cH_2$ and $\cH_\infty$ control problems and derive a graph theoretic condition that guarantees continuous differentiability of the $\cH_\infty$ objective function. In Section~\ref{sec.ls}, we study leader selection in directed networks and, in Section~\ref{sec.cdt}, we design combination drug therapy for HIV treatment. Finally, in Section~\ref{sec:conc}, we conclude the paper and summarize future research directions.

	\vspace*{-2ex}
\section{Problem formulation and background}
	\label{sec:problem:statement}

We first provide background on graph theory and positive systems, describe the system under study, and formulate structured decentralized $\cH_2$ and $\cH_\infty$ optimal control problems.

	\vspace*{-2ex}
\subsection{Background material}
	
\subsubsection*{Notation}
The set of real numbers is denoted by $\R$ and the sets of nonnegative and positive reals are $\R_+$  and $\R_{++}$. The set of $n\times  n$ Metzler matrices (matrices with nonnegative off-diagonal elements) is denoted by $\M^{n \times n}$. We write $A > 0$ ($A \geq 0$) if $A$ has positive (nonnegative) entries and $A \succ 0$ ($A\succcurlyeq0$) if $A$ is symmetric and positive (semi)-definite. We define the sparsity pattern of a vector $u$, $\sparse(u)$, as the set of indices for which $u_i$ is nonzero, $\norm{u}_1 \DefinedAs \sum_i |u_i|$ is the $\ell_1$ norm, and $K^{\dagger}$: $\bbR^{n \times n} \to \bbR^m$ is the adjoint of a linear operator $K$: $\bbR^m \to \bbR^{n \times n}$ if it satisfies,
$
	\inner{X}{K(u)}
	=
	\inner{K^{\dagger}(X)}{u}
$
for all $u \in \bbR^m$ and $X \in \bbR^{n \times n}$.

	\vsp
	
\begin{definition}[Graph associated with a matrix]
$\mc G(A)=(\mc V,\mc E)$ is the graph associated with a matrix $A \in \R^{n\times n}$, with the set of nodes (vertices) $\mc V \DefinedAs \{1,\ldots,n\}$ and the set of edges $\mc E \DefinedAs \{(i,j) | \, A_{ij} \not= 0 \}$, where $(i,j)$ denotes an edge pointing from node $j$ to node $i$.
\end{definition}

	 \vsp
	
\begin{definition}[Strongly connected graph]
A graph $(\mc V,\mc E)$ is strongly connected if there is a directed path between any two distinct nodes in $\mc V$.
\end{definition}

	\vsp
	
\begin{definition}[Weakly connected graph]
A graph $(\mc V,\mc E)$ is weakly connected if replacing its edges with undirected edges results in a strongly connected graph.
\end{definition}

	\vsp
	
	\begin{definition}[Balanced graph]
A graph $(\mc V,\mc E)$ is balanced if, for every node $i \in {\mc V}$, the sum of edge weights on the edges pointing {\em to\/} node $i$ is equal to the sum of edge weights on the edges pointing {\em from\/} node $i$.
	\end{definition}

	\vsp

\begin{definition} A dynamical system is positive if, for any nonnegative initial condition and any nonnegative input, the output is nonnegative for all time. A linear \mbox{time-invariant system,}
\[
	\ba{rcl}
	\dot x
	& \!\!\! = \!\!\! &
	A \, x
	\; + \;
	B \, d
	\\[0.1cm]
	z
	& \!\!\! = \!\!\! &
	C \, x,
	\ea
\]
is positive if and only if $A \in \mathbb{M}^{n \times n}$, $B \geq 0$, and $C \geq 0$.
\end{definition}

	\vsp
	
We now state three lemmas that are useful for the analysis of positive LTI systems.

	\vsp
	
\begin{lemma}[from~\cite{berple94}]
Let $A\in\M^{n \times n}$ and let $Q \in \bbR^{n \times n}$ be a positive definite matrix with nonnegative entries. Then
	\bi
	\item[(a)] $\mre^A\geq0$;
	
	\item[(b)] For Hurwitz $A$, the solution $X$ to the algebraic Lyapunov equation,
$$
	AX
	\;+\;
	XA^T
	\;+\;
	Q
	\;=\;
	0,
$$
is elementwise nonnegative.
	\ei
\label{lem:gramians}
\end{lemma}

	\vsp
		
\begin{lemma}
\label{lem:singvec}
The left and right principal singular vectors, $w$ and $v$, of $A\in\bbR_{+}^{n \times n}$ are nonnegative. If $A \in\bbR_{++}^{n \times n}$, $w$ and $v$ are positive and unique.
\end{lemma}

	\vsp
	
\begin{IEEEproof}
The result follows from the application of the Perron theorem~\cite[Theorem 8.2.11]{horn1985matrix} to $AA^T$ and $A^T A$.
\end{IEEEproof}	

	\vsp
	
	\begin{lemma}[from~\cite{colaneri2014convexity}]
\label{lem:colaneri}
Let $b,c \in \bbR^n$ be nonnegative. Then, for any $t \geq 0$,
\be
c^T  \mre^{(A \, + \, K(u))t}  \,b
\non
\label{eq.cAb}
\ee
is a convex function of $u$.
\end{lemma}
	
	\vspace*{-2ex}
\subsection{Decentralized optimal control} 
	\label{sec.dc}

We consider a class of control problems
	\be
	\label{eq.sys}
	\ba{rcl}
	\dot x
	& \!\!\! = \!\!\! &
	\left(A \, + \, K(u)\right)x
	\; + \;
	B \, d
	\\[0.05cm]
	z
	& \!\!\! = \!\!\! &
	C \, x,
	\ea
	\ee
where $x (t) \in \bbR^n$ is the state vector, $z (t) \in \bbR^m$ is the performance output, $d (t)  \in \bbR^p$ is the disturbance input, and $u \in \bbR^m$ is the control input. Since control enters into the dynamics in a multiplicative fashion, optimal design of $u$ for system~\eqref{eq.sys} is, in general, a challenging nonconvex problem. In what follows, we introduce an assumption which implies that system~\eqref{eq.sys} is positive for any $u$. As we demonstrate in Section~\ref{sec.conv}, under this assumption both the $\cH_2$ and $\cH_\infty$ structured decentralized optimal control problems are convex.

This class of problems can be used to model a variety of control challenges that arise in, e.g., chemical reaction networks and transportation networks. In this paper, we consider leader selection in directed consensus networks as well as combination drug therapy design for HIV treatment.

	\vsp
	
\begin{assumption} \label{ass.pos}
The matrix $A$ in~\eqref{eq.sys} is Metzler, the matrices $B$ and $C$ are nonnegative, and the diagonal matrix $K(u) \DefinedAs \diag \, (Du)$ with $D \in \bbR^{n \times m}$ is a linear function of $u$.
 \end{assumption}

Our objective is to design a stabilizing {\em diagonal\/} matrix $K(u)$ that minimizes amplification from $d$ to $z$. To quantify the average effect of the impulsive disturbance input $d$, we consider the $L_2$ norm of the resulting impulse response, i.e.,
	\begin{subequations} \label{eq.J}
	\be
	J_2(u)
	\DefinedAs 
	\ds\int_{0}^{\infty} \trace 
	\left(
	C \, \mre^{( A + K(u) ) t} B \, B^T \mre^{( A + K(u) )^T t} \, C^T
	\right)
	\mrd t.
	\label{eq.defJ2}
	\ee
This performance metric is equivalent to the square of the $\cH_2$ norm of system~\eqref{eq.sys}~which also has a well-known stochastic interpretation~\cite{dulpag00}. To quantify the worst-case input-output amplification of~\eqref{eq.sys}, we consider the $\mc H_\infty$ norm, defined~as,
	\be
	J_\infty(u)
	\; \DefinedAs \;
	\sup_{\omega \, \in \, \R} 
	~
	\bar\sigma\left(C \, (j\omega I-A-K(u))^{-1}B \right),
	\label{eq.defJinf}
	\ee
	\end{subequations}
where $\bar\sigma (\cdot)$ denotes the largest singular value of a given matrix when $A + K(u)$ is Hurwitz and $\infty$ otherwise. To limit the size of the control input $u$ and promote desired structural properties, we consider the regularized optimal control problem,
\be
\ba{rl}
\minimize\limits_u & J(u) ~+~ g(u) 
	\\[0.1cm]
\subject & A \;+\; K(u) ~~\mbox{Hurwitz}.
\ea
\label{pr.gen}
\ee
The regularization function $g$ in~\eqref{pr.gen} can be any convex function, e.g., a quadratic penalty $u^TRu$ with $R \succ 0$ to limit the magnitude of $u$, an $\ell_1$ penalty to promote sparsity of $u$, or the indicator function associated with a convex set $\cC$ to ensure that $u \in \cC$. We refer the reader to~\cite{linfarjovTAC13admm,jovdhiEJC16} for an overview of recent uses of regularization in control-theoretic problems. 

	\vsp
	
We now review some recent results. Under Assumption~\ref{ass.pos} the matrix $A + K(u)$ is a Metzler and its largest eigenvalue is real and a convex function of $u$~\cite{cohen1981convexity}. Recently, it has been shown that the weighted $\cL_1$ norm of the response of system~\eqref{eq.sys} from a nonnegative initial condition $x_0 \geq 0$,
\[
\ds\int_{0}^T  c^T x(t) \,\mathrm{d}t
\]
is a convex function of $u$ for every $c \in \R^n_+$~\cite{colaneri2014convexity,ranber14}. Furthermore, the approach in~\cite{tanlan11} can be used to cast the problem of {\em unstructured\/} decentralized $\cH_\infty$ control of positive systems as a semidefinite program (SDP) and~\cite{ran16} can be used to cast it as a Linear Program (LP). However, both the SDP and the LP formulations require a change of variables that does not preserve the structure of $K(u)$. Consequently, it is often difficult to design controllers that are feasible for a given noninvertible operator $K$ or to impose structural constraints or penalties on $u$.
	
	\vspace*{-1ex}
\section{Convexity of optimal control problems}
	\label{sec.conv}

We next establish the convexity of the $\cH_2$ and $\cH_\infty$ norms for systems that satisfy Assumption~\ref{ass.pos}, derive a graph theoretic condition that guarantees continuous differentiability of $J_\infty$, and develop a customized algorithm for solving optimization problem~\eqref{pr.gen} even in the absence of differentiability.

	\vspace*{-2ex}
\subsection{Convexity of $J_2$ and $J_\infty$}

We first establish convexity of the $\cH_2$ optimal control problem and provide the expression for the gradient of $J_2$.

	\vspa
	\begin{proposition}
\label{prop.gradJ2}
Let Assumption~\ref{ass.pos} hold and let $\Acl (u) \DefinedAs A + K(u)$ be a Hurwitz matrix. Then, $J_2$ is a convex function of $u$ and its gradient is given by
\be
	\nabla J_2 (u)
	\;=\;
	2 K^{\dagger}(X_cX_o),
\label{eq.grad}
\ee
where $X_c$ and $X_o$ are the controllability and observability gramians of the closed-loop system~\eqref{eq.sys},
\begin{subequations}
\begin{IEEEeqnarray}{rcl}
	\Acl (u) \, X_c
	\;+\;
	X_c \, \Acl^T (u) 
	\;+\;
	B B^T
	& ~ = ~ &
	0
	\label{eq.Xc}
	\\[0.cm]
	\Acl^T (u) \, X_o
	\; + \;
	X_o \, \Acl (u) 
	\;+\;
	C^T  C
	& ~ = ~ &
	0.
	\label{eq.Xo}
\end{IEEEeqnarray}
\end{subequations}
\end{proposition}
	
	\vspa

\begin{IEEEproof}
We first establish convexity of $J_2 (u)$ and then derive its gradient. The square of the $\mc H_2$ norm is given by,
\[
J_2(u)
\;=\;
\left\{
\ba{rl}
\inner{C^T C}{X_c},
&
\Acl(u) \mbox{ Hurwitz}
\\[0.1cm]
\infty,
&
\mbox{otherwise},
\ea
\right.
\]
where the controllability gramian $X_c$ of the closed-loop system is determined by the solution to Lyapunov equation~\eqref{eq.Xc}. For Hurwitz $\Acl (u)$, $X_c$ can be expressed as,
\[
	X_c
	\;=\;
	\ds\int_{0}^{\infty} \mre^{\Acl (u) t} BB^T  \mre^{\Acl^T (u) t} \, \mathrm{d}t.
\]
Substituting into $\inner{C^TC}{X_c}$ and rearranging terms yields,
\[
	\ba{rcl}
	J_2(u)
	& \!\!\! = \!\!\! &
	\ds\int_{0}^\infty\norm{C \, \mre^{\Acl (u) t} B}_F^2 
	\, 
	\mathrm{d}t
	\\[0.3cm]
	& \!\!\! = \!\!\! &
	\ds\int_{0}^\infty \ds\sum_{i, \, j} \left(c_i^T \, \mre^{\Acl (u) t} \, b_j\right)^2
	\mathrm{d} t,
	\ea
\]
where $c_i^T $ is the $i$th row of $C$ and $b_j$ is the $j$th column of $B$.
From Lemma~\ref{lem:colaneri}, it follows that
	$
	c^T  \mre^{\Acl (u) t}  \, b
	$
is a convex function of $u$ for nonnegative vectors $c$ and $b$. Since the range of this function is $\bbR_{+}$ and $(\cdot)^2$ is nondecreasing over $\bbR_{+}$, the composition rules for convex functions~\cite{boyvan04} imply that $(c_i^T  \mre^{\Acl (u) t}  \, b_j)^2$ is convex in $u$. Convexity of $J_2(u)$ follows from the linearity of the sum and integral operators.

	\vsp

To derive $\nabla J_2$, we form the associated Lagrangian,
	\[
	\ba{rcl}
	\mc L(u,X_c,X_o) 
	& \!\!\! = \!\!\! & 
	\inner{C^T C}{ X_c} ~ +
	\\[0.1cm]
	& \!\!\!  \!\!\! & 
	\inner{X_o}{\Acl (u) X_c +X_c\Acl^T (u)  + BB^T},
	\ea
	\]
where $X_o$ is the Lagrange multiplier associated with equality constraint~\eqref{eq.Xc}. Taking variations of $\mc L$ with respect to $X_o$ and $X_c$ yields Lyapunov equations~\eqref{eq.Xc} and~\eqref{eq.Xo} for controllability and observability gramians, respectively. Using $\Acl (u) = A + K(u)$ and the adjoint of $K$, we rewrite the Lagrangian as
	\[
	\ba{rcl}
	\mc L(u,X_c,X_o) 
	& \!\!\! = \!\!\! & 
 	2\inner{K^\dagger (X_cX_o)}{u}  \,+\, \inner{C^T C}{ X_c}
	\; + 
	\\[0.1cm]
	& \!\!\!  \!\!\! & 
	\inner{X_o}{AX_c \, + \, X_cA^T  \, + \, BB^T}.
	\ea
	\]
Taking the variation of $\mc L$ with respect to $u$ yields~\eqref{eq.grad}.
\end{IEEEproof}

	\vspa

\begin{remark}
The quadratic cost
\[
	\int_{0}^T x^T (t) \, C^T C \, x(t) \, \mathrm{d}t,
\]
is also convex over a finite or infinite time horizon for a piecewise constant $u$. This follows from~\cite[Lemma 4]{colaneri2014convexity} and suggests that an approach inspired by the Model Predictive Control (MPC) framework can be used to compute a time-varying optimal control input for a finite horizon problem.
\end{remark}

	\vspa

    \begin{remark}
The expression for $\nabla J_2$ in Proposition~\ref{prop.gradJ2} remains valid for any linear system and any linear operator $K$: $\bbR^m \to \bbR^{n \times n}$. However, convexity of $J_2$ holds under Assumption~\ref{ass.pos} and is not guaranteed in general.
    \end{remark}	
	
We now establish the convexity of the $\cH_\infty$ control problem and provide expression for the subdifferential set of $J_\infty$.

	\vspa
\begin{proposition}
\label{prop.gradJinf}
Let Assumption~\ref{ass.pos} hold and let $\Acl (u) \DefinedAs A + K(u)$ be a Hurwitz matrix. Then, $J_\infty$ is a convex function of $u$ and its subdifferential set is given by
\be	
	\label{eq.grad_hinf}
	\!\!
	\ba{rcl}
	\partial J_\infty(u)
	& \!\!\! = \!\!\! &
	\Big\{ 
	\sum_i \alpha_i \, K^{\dagger}
	\!
	\left(\Acl^{-1} (u) \, B \, v_i \, w_i^T \, C  \Acl^{-1} (u) \right)
	| \; 
	\\[0.1cm]
	& \!\!\! \!\!\! &
	~~
	w_i^T (C \Acl^{-1} (u) B) \, v_i 
	= 
	J_\infty(u), \; \alpha \in \mathcal{P} \Big\},
	\ea
\ee
where $K^{\dagger}$ is the adjoint of the operator $K$ and $\cal P$ is the simplex, ${\cal P} \DefinedAs \{\alpha_i \, | \; \alpha_i \geq 0, \, \sum_i \alpha_i \, = \, 1\}$.
\end{proposition}

	\vspa
	
\begin{IEEEproof}
We first establish convexity of $J_\infty (u)$ and then derive the expression for its subdifferential set. For positive systems, the $\cH_\infty$ norm achieves its largest value at $\omega = 0$~\cite{tanlan11} and from~\eqref{eq.defJinf} we thus have $J_\infty(u) = \bar \sigma(-C\Acl^{-1} (u) B)$. To show convexity of $J_\infty (u)$, we show that $-C\Acl^{-1} (u) B$ is a convex and nonnegative function of $u$, that $\bar\sigma(X)$ is a convex and nondecreasing function of a nonnegative argument $X$, and leverage the composition rules for convex functions~\cite{boyvan04}.

	\vsp

Since $\Acl (u)$ is Hurwitz, its inverse can be expressed as
	\be
	-\Acl^{-1} (u)
	\;=\;
	\ds\int_{0}^{\infty} \mre^{\Acl (u) t} \, \mathrm{d}t.
	\label{eq.Ainv}
	\ee
Convexity of $c_i\mre^{\Acl (u) t} \, b_j$ by Lemma~\ref{lem:colaneri} and linearity of integration implies that each element of the matrix
\[
	-C \, \Acl^{-1} (u) \, B
	\;=\;
	C \ds\int_{0}^{\infty} \mre^{\Acl (u) t} \, \mathrm{d}t \, B
\]
is convex in $u$ and, by part (a) of Lemma~\ref{lem:gramians}, nonnegative.

	\vsp
	
The largest singular value $\bar{\sigma}(X)$ is a convex function of the entries of $X$~\cite{boyvan04},
\be
	\bar{\sigma}(X)
	\;=\;
	\sup_{\norm{w} = 1, \, \norm{v} = 1} w^T Xv,
\label{eq.smax}
\ee
and Lemma~\ref{lem:singvec} implies that the principal singular vectors $v_i$ and $w_i$ that achieve the supremum in~\eqref{eq.smax} are nonnegative for $X \geq 0$. Thus,
\[
	w_i^T
	(X \,+\, \beta \, \mre_j \mre_k^T)
	\,
	v_i
	\;\geq\;
	 w_i^T Xv_i
\]
for any $\beta \geq 0$, thereby implying that $\bar\sigma(X)$ is nondecreasing over $X \geq 0$. Since each element of $-C\Acl^{-1}(u)B \geq 0$ is convex in $u$, these properties of $\bar \sigma (\cdot)$ and the composition rules for convex functions~\cite{boyvan04} imply convexity of $J_\infty(u) = \bar \sigma(-C\Acl^{-1}(u)B)$.

	\vsp

To derive $\partial J_\infty(u)$, we note that the subdifferential set of the supremum over a set of differentiable functions,
\[
    f(x)
    \;=\;
    \sup_{i \, \in \, \mc I}
    ~
    f_i(x),
\]
is the convex hull of the gradients of the functions that achieve the supremum~\cite[Theorem 1.13]{shor2012minimization},
\[
    \partial f(x)
    \;=\;
    \Big\{\ds\sum_{j \,|\, f_j(x) \,=\, f(x)} \alpha_j\nabla f_j(x)
    \, | \;
    \alpha \in {\cal P}
    \Big\}.
\]
Thus, the subgradient of $\bar{\sigma}(X)$ with respect to $X$ is given by
\[
	\ba{c}
	\partial\, \bar{\sigma}(X)
	\;=\;
	\left\{ \sum_j \alpha_j w_jv_j^T  
	\, | \; 
	w_j^T Xv_j 
	\, = \, 
	\bar{\sigma}(X), 
	\;
	\alpha \, \in \, \mathcal{P} 
	\right\}.
	\ea
\]
Finally, the matrix derivative of $X^{-1}$ in conjunction with the chain rule yield~\eqref{eq.grad_hinf}.
\end{IEEEproof}	
	
	\begin{remark}
	For a positive system all induced norms are a convex function of the transfer function matrix evaluated at zero frequncy $-CA_{\text{cl}}^{-1}(u)B$ and Lemma~\ref{lem:colaneri} implies that all induced norms of system~\eqref{eq.sys} are a convex function of $u$. We show convexity of the $\cH_\infty$ norm as it is of particular interest in our study and the proof facilitates the derivation of the gradient.
	\end{remark}
	
	\begin{remark} \label{rem.mono}
	The adjoint of the linear operator $K$, introduced in Assumption~\ref{ass.pos}, with respect to the standard inner product is $K^\dagger(X) = D^T\diag(X)$. For positive systems, Lemma~\ref{lem:gramians} implies that the gramians $X_c$ and $X_o$ are nonnegative matrices. Thus, the diagonal of the matrix $X_cX_o$ is positive and it follows that $J_2$ is a monotone function of the diagonal matrix $K(u)$. Similarly, $-\Acl^{-1}(u)B \, v_i$ and $-w_i^TC\Acl^{-1}(u)$ are nonnegative and thus $J_\infty$ is also a monotone function of $K(u)$.
	\end{remark}

	\vspace*{-2ex}
\subsection{Differentiability of the $\cH_\infty$ norm}

In general, the $\cH_\infty$ norm is a nondifferentiable function of the control input $u$. Even though, under Assumption~\ref{ass.pos}, the decentralized $\cH_\infty$ optimal control problem~\eqref{pr.gen} for positive systems is convex, it is still difficult to solve because of the lack of differentiability of $J_\infty$. Nondifferentiable objective functions often necessitate the use of subgradient methods, which can converge slowly to the optimal solution. 

	\vsp

In what follows, we prove that $J_\infty$ is a continuously differentiable function of $u$ for weakly connected $\cG(A)$. Then, by noting that $J_\infty$ is nondifferentiable only when $\cG(A)$ contains disconnected components, we develop a method for solving~\eqref{pr.gen} that outperforms the standard subgradient algorithm.

	\vspa
	
\subsubsection{Differentiability under weak connectivity}

In this section, we assume that the matrices $B$ and $C$ are square and that their main diagonals are positive. To show the result, we first require two technical lemmas.

	\vspa
	
\begin{lemma} \label{lem.scon}
	Let $M \geq 0$ be a matrix whose main diagonal is strictly positive and whose associated graph $\cG(M)$ is weakly connected. Then, the graphs associated with $\cG(MM^T)$ and $\cG(M^TM)$ have self loops and are strongly connected.
\end{lemma}
	
	\vspa
	
\begin{IEEEproof}
Positivity of the main diagonal of $M$ implies that if $M_{ij}$ is nonzero, then $(M^TM)_{ij}$ and $(MM^T)_{ij}$ are nonzero; by symmetry, $(M^TM)_{ji}$ and $(MM^T)_{ji}$ are also nonzero. Thus, $\cG(M^TM)$ and $\cG(MM^T)$ contain all the edges $(i,j)$ in $\cG(M)$ as well as their reversed counterparts $(j,i)$. Since $\cG(M)$ is weakly connected, $\cG(M^TM)$ and $\cG(MM^T)$ are strongly connected. The presence of self loops follows directly from the positivity of the main diagonal of $M$.
\end{IEEEproof}

	\vspa
	
\begin{lemma} \label{lem.sing}
	Let $M \geq 0$ be a matrix whose main diagonal is strictly positive and whose associated graph $\cG(M)$ is weakly connected. Then, the principal singular value and the principal singular vectors of $M$ are unique.
\end{lemma}

	\vspa

\begin{IEEEproof}
	Note that $\cG(M^k)$ has an edge from $i$ to $j$ if $M$ contains a directed path of length $k$ from $i$ to $j$~\cite[Lemma 1.32]{bullo2009distributed}. Since $\cG(MM^T)$ and $\cG(M^TM)$ are strongly connected with self loops, Lemma~\ref{lem.scon} implies the existence of $\bar k$ such that $(M^TM)^{k} > 0$ and $(MM^T)^{k} > 0$ for all $k \geq \bar k$, and Perron Theorem~\cite[Theorem 8.2.11]{horn1985matrix} implies that $(M^TM)^k$ and $(MM^T)^k$ have unique maximum eigenvalues \mbox{for all $k \geq \bar k$.}
	
		\vspa
	
	The eigenvalues of $(M^TM)^k$ and $(MM^T)^k$ are related to the singular values of $M$ by, 
	\[
	\lambda_i((M^TM)^k) 
	\, = \,
	\lambda_i((MM^T)^k)
	\, = \, 
	( \sigma_i(M) )^{2k},
	\]
and the eigenvectors of $(M^TM)^k$ and $(MM^T)^k$ are determined by the right and the left singular vectors of $M$, respectively. Since the principal eigenvalue and eigenvectors of these matrices are unique, the principal singular value and the associated singular vectors of $M$ are also unique.
\end{IEEEproof}

	\vspa
	
\begin{theorem} \label{thm.weak}
Let Assumption~\ref{ass.pos} hold, let $\Acl (u) \DefinedAs A + K(u)$ be a Hurwitz matrix, and let matrices $B$ and $C$ have strictly positive main diagonals. If the graph $\cG(A)$ associated with $A$ is weakly connected, $J_\infty(u)$ is continuously differentiable.
\end{theorem}

	\vspa
	
\begin{IEEEproof}
Lemma~\ref{lem:gramians} implies that~$\mre^{\Acl (u)} \geq 0$. From~\cite[Lemma 1.32]{bullo2009distributed}\label{lem.bullo}, $\cG(M^k)$ has an edge from $i$ to $j$ if there is a directed path of length $k$ from $i$ to $j$ in $\cG(M)$. 
Weak connectivity of $\cG(A)$ implies weak connectivity of $\cG(\tilde A)$, ${\cG}(\tilde{A}^k)$, $\mre^{\Acl (u) t}$ and, by~\eqref{eq.Ainv}, of ${\cG}(-\Acl^{-1} (u))$.

	\vsp

Since $\Acl (u)$ is Hurwitz and Metzler, its main diagonal must be strictly negative; otherwise, $\tfrac{\mathrm{d}}{\mathrm{d}t}x_i \geq 0$ for some $x_i$, contradicting stability and thus the Hurwitz assumption. Equation~\eqref{eq.Ainv} and Lemma~\ref{lem:gramians} imply $\Acl^{-1}(u) \leq 0$ and, since $\Acl(u)$ is Metzler, $\Acl^{-1}(u)\Acl(u) = I$ can only hold if the main diagonal of $-\Acl^{-1}(u)$ is strictly positive.

	\vsp
	
Moreover, since the diagonals of $B$ and $C$ are strictly positive, $\cG(-C\Acl^{-1} (u) B)$ is weakly connected and the diagonal of $-C\Acl^{-1} (u)B$ is also strictly positive. Thus, Lemma~\ref{lem.sing} implies that the principal singular value and singular vectors of $-C\Acl^{-1} (u) B$ are unique, that~\eqref{eq.grad_hinf} is unique for each stabilizing $u$, and that $J_{\infty} (u)$ is continuously differentiable.
\end{IEEEproof}

	\vspa

\subsubsection{Nondifferentiability for disconnected $\cG(A)$}

Theorem~\ref{thm.weak} implies that under a mild assumption on $B$ and $C$, $J_\infty$ is only nondifferentiable when the graph associated with $A$ has disjoint components. Proximal methods and its accelerated variants~\cite{becteb09} generalize gradient descent to nonsmooth problems when the proximal operator of the nondifferentiable term in the objective function is readily available. However, since there is no explicit expression for the proximal operator of $J_\infty$,  in general we have to use subgradient methods to solve~\eqref{pr.gen}.

	\vsp

To a large extent, subgradient methods are inefficient because they do not guarantee descent of the objective function. However, under the following mild assumption, the subgradient of $J_\infty$, $\partial J_\infty$, can be conveniently expressed and a descent direction can be obtained by solving a linear program. 
 	
	\vsp
	
\begin{assumption}\label{ass.weak}
Without loss of generality, let $\Acl(u)$ be permuted such that $\Acl(u) = \blkdiag \, (\Acl^1(u),\dots,\Acl^m(u))$ is block diagonal and let $\cG(\Acl^i(u))$ be weakly connected for every $i$. Moreover, the matrices $B = \blkdiag \, (B^1, \dots, B^m)$ and $C = \blkdiag \, (C^1, \dots, C^m)$ are block diagonal and partitioned conformably with the matrix $\Acl(u)$.
\end{assumption}

	\vspa
	
	\begin{theorem} \label{thm.blk}
Let Assumptions~\ref{ass.pos} and~\ref{ass.weak} hold and let $\Acl(u) \DefinedAs A + K(u)$ be a Hurwitz matrix.  Then,
\begin{subequations} \label{eq.blk}
\be
	J_\infty(u)
	\;=\;
	\max_i
	\,
	J^i_\infty(u)	
	\label{eq.blkJ},
\ee
where $J_\infty^i(u) \DefinedAs \bar\sigma (C^i(\Acl^i(u))^{-1}B^i)$. Moreover, every element of the subgradient of $J_\infty(u)$ can be expressed as the convex combination of a finite number of vectors $f^j \DefinedAs \nabla J_\infty^j(u)$ corresponding to the gradients of the functions $J_\infty^j(u)$ that achieve the maximum in~\eqref{eq.blkJ}, i.e., $J_\infty(u) = J_\infty^j(u)$,
\be	
	\partial J_\infty(u)
	\; = \;
	\left\{
	F\alpha
	~|~ \alpha \in \mathcal{P} \right\}
	\label{eq.blknablaJ}
\ee
\end{subequations}
where the columns of $F$ are given by $f^j$ and $\cal P$ is the simplex.
\end{theorem}
		
\begin{IEEEproof}
	Since $\Acl(u)$ is a block diagonal matrix, so is $\Acl^{-1}(u)$ and Assumption~\ref{ass.weak} implies that $-C\Acl^{-1} (u) B = \blkdiag \, (-C^i(\Acl^i(u))^{-1}B^i)$ is also block diagonal. Thus,
	\[
		J_\infty(u)
		\;=\;
		\bar\sigma(-C\Acl^{-1}(u)B)
		\;=\;
		\ds\max_i \; \bar\sigma(-C^i(\Acl^i(u))^{-1}B^i),
	\]
which proves~\eqref{eq.blkJ}. Theorem~\ref{thm.weak} implies that each $J^i_\infty(u)$ is continuously differentiable which establishes~\eqref{eq.blknablaJ}.
\end{IEEEproof}

	\vspa

When $g$ is differentiable, we leverage the above convenient expression for $\partial J_\infty$ to select an element of the subdifferential set which, with an abuse of terminology,  we call the {\em optimal\/} subgradient. The optimal subgradient is guaranteed to be a descent direction for~\eqref{pr.gen} and it is defined as the member of $\partial (J_\infty(u) + g(u))$ that solves
	\begin{subequations}
	\label{eq.optsg}
	\begin{IEEEeqnarray}{rl}
	\minimize_{v,\,\alpha}
	~&
	\max_j\; (v^T(f^j \,+\, \nabla g (u)))
	\label{eq.gs1}
	\\[0.1cm]
	\subject
	~&
	v
	\,=\,
	-(F\alpha
	\,+\,
	\nabla g(u)),
	~~
	\label{eq.gs2}
	\alpha \,\in\, {\cal P}
	\\[0.1cm]
	~&
	v^T(f^j \,+\, \nabla g (u))
	\,<\,
	0,
	~~
	\mbox{for~all}~j,
	\label{eq.gs4}
	\end{IEEEeqnarray}
	\label{eq.gs}
	\end{subequations}
	
\vspace{-3ex}
\noindent
where $F$ and $f^j$ are defined as in Theorem~\ref{thm.blk}. By~\eqref{eq.blkJ}, $J_\infty(u)$ is the maximum of differentiable functions $J_\infty^i(u)$ and problem~\eqref{eq.gs} forms a search direction using the gradients of the functions $J_\infty^j$ that achieve that maximum. While constraint~\eqref{eq.gs2} ensures that $v \in \partial J_\infty(u) + \nabla g(u)$,~\eqref{eq.gs4} ensures that $v$ is a descent direction for each $J_\infty^j$ and thereby guarantees that $v$ is a descent direction for $J_\infty$. Finally, objective function~\eqref{eq.gs1} is the maximum of the directional derivatives of $J_\infty^j$ in the direction $v$, i.e., the directional derivative of the objective function in~\eqref{pr.gen} in the search direction.

	\vsp
	
Problem~\eqref{eq.gs} can be solved efficiently because it is a linear program. Moreover, the optimality condition for~\eqref{pr.gen}, $\partial J_\infty(u) + \nabla g(u) \ni 0$, can be checked by solving a linear program to verify the existence of an $\alpha \in {\cal P}$ such that $F\alpha + \nabla g(u) = 0$.

	\vsp

	\subsubsection{Customized algorithm}
Ensuring a descent direction enables principled rules for step-size selection and makes problem~\eqref{pr.gen} with nondifferentiable $g$ tractable via the augmented-Lagrangian-based approaches. Reformulation of~\eqref{pr.gen},
\be
	\label{pr.split}
	\ba{rl}
	\minimize\limits_{u, \, v}
	&
	J_\infty(u)
	\;+\;
	g(v)
	\\[0.1cm]
	\subject
	&
	u 
	\, - \, v
	\;=\;
	0,
	\ea
\ee
leads to the associated augmented Lagrangian
\be
	\label{eq.al}
	\non
	\cL_\mu(u,v;\lambda)
	\DefinedAs
	J_\infty(u) \,+\, g(v) \,+\, \lambda^T(u - v) \,+\, \tfrac{1}{2\mu} \, \norm{u - v}^2,
	\notag
\ee
where $v$ is an auxiliary variable and $\mu$ is a positive parameter. Formulation~\eqref{pr.split} is convenient for the alternating direction method of multipliers (ADMM)~\cite{boyparchupeleck11}, which minimizes $\cL_\mu$ separately over $u$ and $v$  and updates $\lambda$ until convergence. ADMM is highly sensitive to the choice of $\mu$ and it may require many iterations to converge. In contrast, the more mature and robust method of multipliers (MM)~\cite{ber99} has effective rules for adaptively updating $\mu$ which leads to faster convergence. It is difficult to directly apply MM to~\eqref{pr.split} because it requires joint minimization of $\cL_\mu$ over $(u,v)$ and both $g$ and $J_\infty$ are nondifferentiable. However, when the proximal operator of $g$ is readily available, e.g., when $g$ is the $\ell_1$ norm or an indicator function of a convex set with simple projection~\cite{parboy13}, explicit minimization over $v$ is achieved by $v^\star_\mu (u,\lambda) = \prox_{\mu g}(u + \mu\lambda)$. Substitution of $v^\star_\mu (u,\lambda)$ into $\cL_\mu$ yields the proximal augmented Lagrangian~\cite{dhikhojovTAC16},
  \be
    \cL_\mu(u, v^\star_\mu (u,\lambda); \lambda)
	\; = \;
    J_\infty(u)
    \; + \;
    M_{\mu g}
    (u \, + \, \mu \lambda)
    \; - \;
    \tfrac{\mu}{2} \, \| \lambda \|^2,
    \label{eq.alprox}
    \notag
    \ee
where $M_{\mu g}$ is the Moreau envelope of $g$ and is continuously differentiable, even when $g$ is not~\cite{parboy13}; see~\cite{dhikhojovTAC16} for details. Since $J_\infty$ is the only nondifferentiable component of the proximal augmented Lagrangian, the optimal subgradient~\eqref{eq.optsg} can be used to minimize it over $u$. This equivalently minimizes $\cL_\mu (u,v;\lambda)$ over $(u,v)$ and leads to a tractable MM algorithm,
	\begin{eqnarray}
	u^{k+1}
	& \!\! = \!\! &
	\argmin\limits_{u}
	\,
	\cL_{\mu} (u, v^\star_{\mu^k} (u,\lambda^k); \lambda^k)
	\label{eq.MMa1}
	\notag
	\\
	\lambda^{k+1}
	& \!\! = \!\! &
	\lambda^{k}
	\; + \;
		\tfrac{1}{\mu^k}
	 \,
	 (
	 u^{k+1}
	 \,-\,
	 \prox_{\mu^k g}(u^{k+1} + \mu^k \lambda^k)
	 ).
	\label{eq.MMa2}
	\notag
	\end{eqnarray}
This algorithm minimizes the proximal augmented Lagrangian over $u$, updates $\lambda$ via gradient ascent, and it represents an appealing alternative to ADMM for problems of the form~\eqref{pr.gen}. In particular, an adaptive selection of the parameter $\mu$ leads to improved practical performance relative to ADMM~\cite{dhikhojovTAC16}. 

	\vspace*{-2ex}
\section{Leader selection in directed networks}
\label{sec.ls}

We now consider the special case of system~\eqref{eq.sys}, in which the matrix $A$ is given by a graph Laplacian, and study the leader selection problem for directed consensus networks. The question of how to optimally assign a predetermined number of nodes to act as leaders in a network of dynamical systems with a given topology has recently emerged as a useful proxy for identifying important nodes in a network~\cite{patbam10,fitch2013information,fitch2014joint,linfarjovTAC14leaderselection,clabuspoo14,clark2016submodularity}. Even though significant theoretical and algorithmic advances for undirected networks have been made, the leader selection problem in directed networks remains open. 

	\vspace*{-2ex}
\subsection{Problem formulation}

We describe consensus dynamics and state the problem.

	\vsp
	
\subsubsection{Consensus dynamics}

The weighted directed network $\cG(L)$ with $n$ nodes and the graph Laplacian $L$ obeys consensus dynamics in which each node $i$ updates its state $x_i$ using relative information exchange with its neighbors, 
	$$
	\dot x_i 
	\; = \; 
	- \sum_{j \, \in \, {\mc N}_i} \! L_{ij} (x_i \, - \, x_j) 
	\; + \; 
	d_i.
	$$
Here, ${\mc N_i} \DefinedAs \{j \, | \, (i,j) \in {\cal E}\}$, $L_{ij} \geq 0$ is a weight that quantifies the importance of the edge from node $j$ to node $i$, $d_i$ is a disturbance, and the aggregate dynamics are~\cite{xiaboykim07},
\[
	\dot x
	\;=\;
	-Lx
	\;+\;
	d,
\]
where $L$ is the graph Laplacian of the directed network~\cite{cvetkovic1980spectra}.

	\vsp

The graph Laplacian always has an eigenvalue at zero that corresponds to a right eigenvector of all ones, $L\one = 0$. If this eigenvalue is simple, all node values $x_i$ converge to a constant $\bar x$ in the absence of an external input $d$.  When $\cG(L)$ is balanced, $\bar x = (1/n) \, \one^Tx(0)$ is the average of the initial node values. In general,
$
	\bar x
	=
	w^Tx(0),
$
where $w$ is the left eigenvector of $L$ corresponding to zero eigenvalue, $w^TL = 0$. If $\cG(L)$ is not strongly connected, $L$ may have additional eigenvalues at zero and the node values converge to distinct groups whose number is equal to or smaller than the multiplicity of the zero eigenvalue.

	\vsp

\subsubsection{Leader selection}

In consensus networks, the dynamics are governed by relative information exchange and the node values converge to the network average. In the leader selection paradigm~\cite{linfarjovTAC14leaderselection}, certain ``leader'' nodes are additionally equipped with {\em absolute\/} information which introduces negative feedback on the states of these nodes. If suitable leader nodes are present, the dynamical generator becomes a Hurwitz matrix and the states of all nodes asymptotically converge to zero.

	\vsp

The node dynamics in a network with leaders is
	$$
	\dot x_i  
	\; = \;  
	- \sum_{j \, \in \, {\mc N}_i} \! 
	L_{ij} 
	(x_i \, - \, x_j) 
	\; - \; 
	u_i \, x_i 
	\;+\; d_i,
	$$
where $u_i \geq 0$ is the weight that node $i$ places on its absolute information. The node $i$ is a leader if $u_i > 0$, otherwise it is a follower. The aggregate dynamics can be written as
\[
\dot x
	\;=\;
	-(L
	\,+\,
	\diag \, (u)
	)
	\,
	x
	\;+\;
	d,
\]
and placed in the form~\eqref{eq.sys} by taking $A = -L$, $B = C = I$, and $K(u) = -\diag \, (u)$. We evaluate the performance of a leader vector $u \in \bbR^n$ using the $\cH_2$ or $\cH_\infty$ performance metrics $J_2$ or $J_\infty$, respectively. We note that this system is marginally stable in the absence of leaders and much work on consensus networks focuses on driving the {\em deviations\/} from the average node values to zero~\cite{bamjovmitpat12}. Instead, we here focus on driving the node values themselves to zero.

	\vsp

We formulate the combinatorial problem of selecting $N$ leaders to optimize either $\cH_2$ or $\cH_\infty$ norm as follows. 

	\vsp

\begin{problem}
\label{pr.comb}
Given a network with a graph Laplacian $L$ and a fixed leader weight $\kappa$, find the optimal set of $N$ leaders that solves
\[
	\ba{rl}
	\minimize\limits_u
	&
	J(u)
	\\[0.1cm]
	\subject
	&
	\one^T u
	\; = \;
	N\,\kappa,
	~~
	u_i
	\, \in \,
	\{0,\kappa\},
	\ea
\]
where $J$ is a performance metrics described in Section~\ref{sec.dc}, with $A = - L$, $B = C = I$, and $K(u) = \diag \, (u)$.
\end{problem}

	\vsp

In~\cite{fitch2013information,fitch2014joint}, the authors derive explicit expressions for leaders in undirected networks. However, these expressions are efficient only for very few or very many leaders. Instead, we follow~\cite{linfarjovTAC14leaderselection} and develop an algorithm which relaxes the integer constraint to obtain a lower bound on Problem~\ref{pr.comb} and use greedy heuristics to obtain an upper bound.

	\vsp

Considering leader selection in directed networks adds the challenge of ensuring stability. At the same time, we can leverage existing results on leader selection in undirected networks to derive efficient upper bounds on Problem~\ref{pr.comb}.

	\vspace*{-2ex}
\subsection{Stability for directed networks}

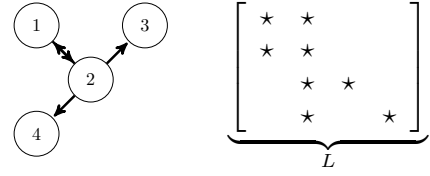
\begin{figure}
\[
	\begin{aligned}
	\raisebox{-11mm}{
	\resizebox{.25\hsize}{!}{
	\input{smallnetwork}}}
	&
	~~~
	&
	\resizebox{.3\hsize}{!}{$
	\underbrace{
         \matbegin \begin{array}{cccc}
                \star & \star &  & \\
                 \star & \star &  & \\
                  & \star & \star & \\
                 & \star &  &  \star
                \end{array}
                \matend
                }_{L}$}
        \end{aligned}
        \]
\caption{A directed network and the sparsity pattern of the corresponding graph Laplacian. This network is stabilized if and only if either node $1$ or node $2$ are made leaders.}
\label{fig.unstable}
\end{figure}

For a vector of leader weights $u$ to be feasible for Problems~\ref{pr.comb}, it must stabilize system~\eqref{eq.sys}; i.e., $-(L + \diag \, (u))$ must be a Hurwitz matrix. When $\cG(L)$ is undirected and connected, any leader will stabilize~\eqref{eq.sys}. However, this is not the case for directed networks. For example, making node $1$ or $2$ a leader stabilizes the network in Fig.~\ref{fig.unstable}, but making nodes $3$ or $4$ a leader does not. Theorem~\ref{thm.stab} provides a necessary and sufficient condition for stability.

	\vspa

\begin{theorem}
\label{thm.stab}
Let $L$ be a weighted directed graph Laplacian and let $u \geq 0$. The matrix $- (L + \diag \, (u))$ is Hurwitz if and only if $w \circ u \not= 0$ for all nonzero $w$ with $w^TL = 0$, where $\circ$ is the elementwise product.
\end{theorem}

	\vspa
		
\begin{IEEEproof}
$(\Leftarrow)$
If $w \circ u = 0$, $w^T\diag \, (u) = 0$. If, in addition, $w^TL = 0$, we have
\be
	-w^T(L \,+\, \diag \, (u)) 
	\;=\;
	 0
	 \label{eq.Lzero},
\ee
and therefore zero is an eigenvalue of $-(L + \diag \, (u))$.

	\vsp
	
$(\Rightarrow)$
Since the graph Laplacian $L$ is row stochastic and $\diag \, (u)$ is diagonal and nonnegative, the Gershgorin circle theorem~\cite{horn1985matrix} implies that the eigenvalues of $-(L + \diag \, (u))$ are at most $0$. To show that $-(L + \diag \, (u))$ is Hurwitz, we show that it has no eigenvalue at zero. Assume there exists a nonzero $w$ such that~\eqref{eq.Lzero} holds. This implies that either
$
    w^TL
    =
    w^T\diag \, (u)
    =
     0
$
 or that
$
    w^TL
    =
    -w^T\diag \, (u).
$
The first case is not possible because, by assumption, $w^T\diag \, (u) = (w \,\circ\, u)^T \not= 0$ for any $w$ such that $w^TL = 0$. If the second case is true, then $w^TLv = -w^T\diag \, (u)v$ must also hold for all $v$. However, if we take $v = \one$, then $w^TL\one = 0$ but $-w^T\diag \, (u)\one$ is nonzero. This completes the proof.
\end{IEEEproof}

	\vspa

	\begin{remark}
Only the set of leader nodes is relevant to the question of stability. If $u$ does not stabilize~\eqref{eq.sys}, no positive weighting of the vector of leader nodes, $\alpha \circ u$ with $\alpha \in \bbR_{++}^N$, will stabilize~\eqref{eq.sys}. Similarly if $u$ stabilizes~\eqref{eq.sys}, every $\alpha \circ u$ will.
\end{remark}

	\vspa
		
	\begin{corollary} \label{cor.str}
If $\cG(L)$ is strongly connected, any choice of leader node will stabilize~\eqref{eq.sys}.
\end{corollary}

	\vspa
	
\begin{IEEEproof}
Since the graph Laplacian associated with a strongly connected graph is irreducible, the Perron-Frobenius theorem~\cite{horn1985matrix} implies that the left eigenvector associated with $-L$ is positive. Thus, $w \circ u \not = 0$ for any nonzero $u$ and system~\eqref{eq.sys} is stable by Theorem~\ref{thm.stab}.
\end{IEEEproof}

\vspa

\begin{remark} \label{rem.Sj}
The condition in Theorem~\ref{thm.stab} requires that there is a path from the set of leader nodes to every node in the network. This can be enforced by extracting disjoint ``leader subsets'' $\cS_j$ which are not influenced by the rest of the network, i.e., $(v^j)^TL = 0$ where $v^j_i = 1$ if $i \in \cS_j$ and $v^j_i = 0$ otherwise, and which are each strongly connected components of the original network. Stability is guaranteed if there is at least one leader node in each such subset $\cS_j$; e.g., for the network in Fig.~\ref{fig.unstable}, there is one leader subset $\cS_1 = \{1,2\}$. By Corollary~\ref{cor.str}, $\cS_1$ contains all nodes when the network is strongly connected.
\end{remark}

	\vspace*{-2ex}
\subsection{Bounds for Problem~\ref{pr.comb}} \label{sec.bds}

To approach combinatorial Problem~\ref{pr.comb}, we derive bounds on its optimal objective value. These bounds can also be used to implement a branch-and-bound approach~\cite{lawwoo66}.

	\vsp
	
\subsubsection{Lower bound}
By relaxing the combinatorial constraint in Problem~\ref{pr.comb}, we formulate the optimization problem,
\be
\label{eq.opt.prob.lb}
	\ba{rl}
	\minimize\limits_u
	&
	J(u)
	\\[0.1cm]
	\subject
	&
	u
	\, \in \,
	\kappa
	\,
	{\cal P}_N,
	\ea
\ee
where $\kappa\,{\cal P}_N \DefinedAs \{ u \, | \, \sum_i u_i = N\kappa,\, u_i \leq \kappa \}$ is the ``capped'' simplex. The results of Section~\ref{sec.conv}, establish the convexity of problem~\eqref{eq.opt.prob.lb}.  Using a recent result on efficient projection onto ${\cal P}_N$~\cite{wancan15}, this problem can be solved efficiently via proximal gradient methods~\cite{becteb09} to provide a lower bound on Problem~\ref{pr.comb}.

	\vsp

When $\cG(L)$ is not strongly connected, additional constraints can be added to enforce the condition in Theorem~\ref{thm.stab} and thus guarantee stability. Let the sets $\cS_j$ denote ``leader subsets'' from which a leader must be chosen, as discussed in Remark~\ref{rem.Sj}. Then, the convex problem
\be
	\label{pr.rel}
	\ba{rl}
	\minimize
	&
	J(u)
	\\[0.1cm]
	\subject
	&
		u
	\, \in \,
	\kappa
	\,
	{\cal P}_N,
	~
	\ds\sum_{i \,\in\,{\cal S}_{j}} u_i
	\, \geq \,
	\kappa,
	~
	\mbox{for~all}~j
	\ea
\ee
relaxes the combinatorial constraint and guarantees stability. 
We denote the resulting lower bounds on the optimal values of the $\cH_2$ and $\cH_\infty$ versions of Problem~\ref{pr.comb} with $N$ leaders by $\Jlb_2(N)$ and $\Jlb_\infty(N)$, respectively.

	\vsp

\subsubsection{Upper bounds for Problem~\ref{pr.comb}} \label{sec.ub}

If $k$ denotes the number of subsets $\cS_j$, a stabilizing candidate solution to Problem~\ref{pr.comb} can be obtained by ``rounding'' the solution to~\eqref{pr.rel} by taking $N$ leaders to contain the largest element from each subset $\cS_j$ and $N - k$ largest remaining elements. The greedy swapping algorithm proposed in~\cite{linfarjovTAC14leaderselection} can further tighten \mbox{this upper bound.} 

	\vsp

Recent work on leader selection in undirected networks can also provide upper bounds for Problem~\ref{pr.comb} when $\cG(L)$ is balanced. The symmetric component of the Laplacian of a balanced graph, $L_s \DefinedAs \frac{1}{2}(L + L^T)$, is the Laplacian of an undirected network. The exact optimal leader set for an undirected network can be efficiently computed when $N$ is either small or large~\cite{fitch2013information,fitch2014joint}. Since the performance of the symmetric component of a system provides an upper bound on the performance of the original system, these sets of leaders will have better performance with $L$ than with $L_s$ for both the $\Ht$~\cite[Corollary 3]{dhijovACC15} and $\Hinf$ norms~\cite[Proposition 4]{dhiwujovNOL16}.

	\vsp
	
Even when $L$ does not represent a balanced network, $J_2$ and $J_\infty$ are respectively upper bounded by the trace and the maximum eigenvalue of $\tfrac{1}{2}(L_s + K)^{-1}$. For small numbers of leaders, they can be efficiently computed using rank-one inversion updates. {A} similar approach was used in~\cite{fitch2013information,fitch2014joint} to derive optimal leaders for undirected networks. Moreover, this approach yields {a} stabilizing set of \mbox{leaders~\cite[Lemma 1]{dhijovACC15}.} 

	\vspace*{-2ex}
\subsection{Additional comments}

We now provide additional discussion on interesting aspects of Problem~\ref{pr.comb}. We first consider the gradients of $J_2$ and $J_\infty$.

\vsp
	
\begin{remark}
When $K(u) = -\diag \, (u)$, we have $\nabla J_2 = -2\,\diag \, (X_cX_o)$. The matrix $X_cX_o$ often appears in model reduction and $(\nabla J_2(u))_i$ corresponds to the inner product between the $i$th columns of $X_c$ and $X_o$. 
\end{remark}
	\vsp
	
\begin{remark}
When $K(u) = -\diag(u)$, $(\partial J_\infty(u))_i$ is given by the product of $-\mre_i^T \Acl^{-1}(u)v$ and $w^T \Acl^{-1}(u)\mre_i$. The former quantifies how much the forcing which causes the largest overall response of system~\eqref{eq.sys} affects node $i$, and the latter captures how much the forcing at node $i$ affects the direction of the largest output response.
\end{remark}

	\vsp
	
The optimal leader sets for balanced graphs are interesting because they are invariant under reversal of all edge directions.

	\vspa
	
\begin{proposition}
	\label{th.cycles}
Let $\cG(L)$ be balanced, let $\hat L \DefinedAs L^T$ so that $\cG(\hat L)$ contains the reversed edges of the graph $\cG(L)$, and let $\hat J_2$ and $\hat J_\infty$ denote the performance metrics~\eqref{eq.J} with $A = -\hat L$, $K(u) = - \diag \, (u)$, and $B = C = I$ as in Problem~\ref{pr.comb}. Then, $J_2(u)= \hat J_2(u)$ and $J_\infty(u) = \hat J_\infty(u)$.
\end{proposition}

	\vspa
	
\begin{IEEEproof}
The controllability gramian of~\eqref{eq.sys} defined with $\Acl = -(L + \diag(u))$ solves Lyapunov equation~\eqref{eq.Xc}, $-(L + \diag(u))X_c - X_c(L + \diag(u))^T + I = 0$, and is also the observability gramian $\hat X_o$ of~\eqref{eq.sys} defined with $\Acl = - (\hat L + \diag(u)) = - (L^T + \diag(u))$ that solves~\eqref{eq.Xo}.  By definition of the $\cH_2$ norm, $\hat J_2(u) = \trace(\hat X_o) = \trace(X_c) = J_2(u)$. Since $\bar \sigma(M) = \bar \sigma(M^T)$, $\hat J_\infty(u) = \bar \sigma(-(\hat L + \diag(u))^{-1}) = \bar \sigma(-(L + \diag(u))^{-1}) = J_\infty(u)$.
\end{IEEEproof}

	\vspa
	
This invariance is intriguing because the space of balanced graphs is spanned by cycles. In~\cite{zelschall13}, the authors explored how undirected cycles affect undirected consensus networks. Proposition~\ref{th.cycles} suggests that directed cycles also play a fundamental role in directed consensus networks.

	\vspace*{-2ex}
\subsection{Computational experiments}
	\label{sec.ex}

Here, we illustrate our approach to Problem~\ref{pr.comb} with $N$ leaders and the weight $\kappa = 1$. The ``rounding" approach that we employ is described in Section~\ref{sec.ub}.

	\vsp

\subsubsection{Synthetic example}

\begin{figure}[t]
\begin{center}
\vspace{-.5cm}
	\begin{tabular}{cc}
	    \!\!\!\!\!
    \!\!\!\!\!
	\subfloat[Balanced Network with $8$ nodes, $12$ edges.]
{
\label{fig.dirnetfig}
	\begin{tabular}{c}
    \!\!\!\!\!     \!\!\!\!\!
        \input{dirnet_bal_med}
    \!\!\!\!\!
    \!\!\!\!\!
    \end{tabular}
    }    
    &
	\subfloat[Performance of optimal leaders and two leader selection techniques.]
	{
    \label{fig.dirnetperf}
  \begin{tabular}{rc}
    \!\!\!\!\!
        \rotatebox{90}{\footnotesize \quad\; $100*(J_2(u)/\Jlb_2(N)-1)$}
    \!\!\!\!\!
    &
        \!\!\!\!\!
    \!\!\!\!\!
    {	\includegraphics[width= .6\columnwidth]{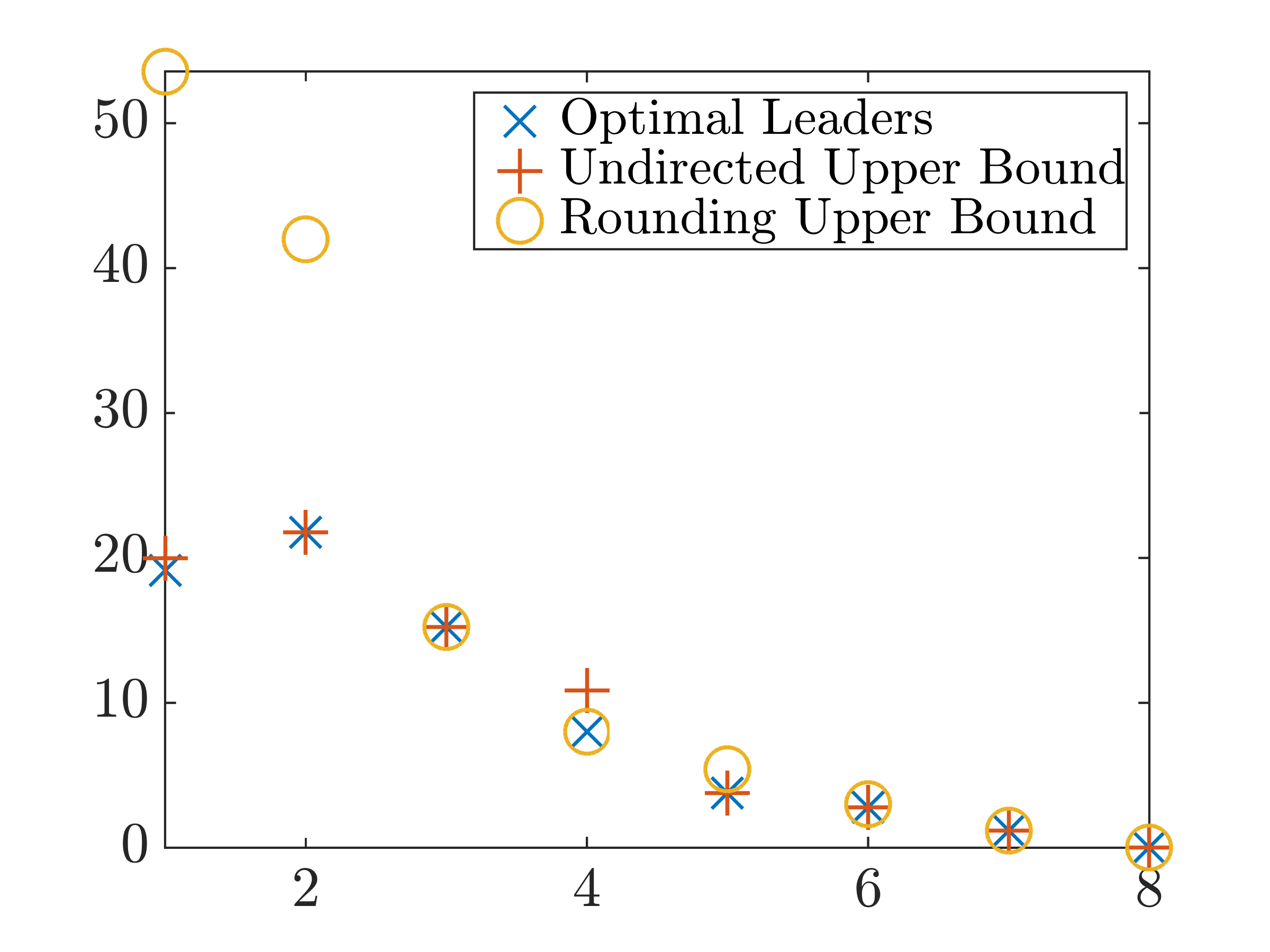}}
    \!\!\!\!\!    \!\!\!\!\!
    \\
    \!\!\!\!\!
        \!\!\!\!\!
    \!\!\!\!\!
    &
    \!\!\!\!\!
    {\footnotesize Number of leaders $N$}
    \!\!\!\!\!
    \end{tabular}}
        \!\!\!\!\!
    \!\!\!\!\!
	\end{tabular}
\end{center}
    \caption{$\cH_2$ performance of optimal leader set (\tc{matlabblue}{blue $\times$}) and upper bounds resulting from ``rounding'' (\tc{matlaborange}{yellow $\circ$}) and the optimal leaders for the undirected network (\tc{matlabred}{red $+$}). Performance is shown as a percent increase in $J_2$ relative to $\Jlb_2(N)$.
    }
	\label{fig.dirnet}
\end{figure}

For the directed network in Fig.~\ref{fig.dirnetfig}, let the edges from node $2$ to node $7$ and from node $7$ to node $8$ have an edge weight of $2$ and let all other edges have unit edge weights. We compare the optimal set of leaders, determined by exhaustive search, to the set of leaders obtained by (i) ``rounding'' the solution to relaxed problem~\eqref{pr.rel}; and by (ii) the optimal selection for the undirected version of the graph via~\cite{fitch2013information,fitch2014joint}, as discussed in Section~\ref{sec.ub}. In Fig.~\ref{fig.dirnetfig}, \tc{matlabblue}{blue node $7$} represents the optimal single leader, \tc{matlaborange}{yellow node $4$} represents the single leader selected by ``rounding", and  \tc{matlabred}{red node $8$} represents the optimal single leader for the undirected network.
In Fig.~\ref{fig.dirnetperf} we show the $\Ht$ performance for $1$ to $8$ leader nodes resulting from different methods. Since in general, we do not know the optimal performance {\em a priori\/}, we plot performance degradation (in percents) relative to the lower bound on Problem~\ref{pr.comb} obtained by solving problem~\eqref{pr.rel}. 

	\vsp
		
Figure~\ref{fig.dirnetperf} shows that neither ``rounding'' (\tc{matlaborange}{yellow $\circ$}) nor the optimal selection for undirected networks (\tc{matlabred}{red $+$}) achieve unilaterally better $\cH_2$ performance (performance of the optimal leader sets are shown in \tc{matlabblue}{blue $\times$}). While the procedure for the undirected networks selects better sets of $1$, $2$, and $5$ leaders relative to ``rounding'', identifying them is expensive except for large or small number of leaders~\cite{fitch2013information,fitch2014joint} and ``rounding'' identifies a better set of $4$ leaders. This suggests that, when possible, {\em both\/} sets of leaders should be computed and the one that achieves better performance should be selected.

	\vsp
	
\subsubsection{Neural network of the worm C.\ Elegans}

We now consider the network of neurons in the brain of the worm C.\ Elegans with $297$ nodes and $2359$ weighted directed edges. The data was compiled by~\cite{watstr98} from~\cite{whisouthobre86}. Inspired by the use of leader selection as a proxy for identifying important nodes in a network~\cite{patbam10,fitch2013information,fitch2014joint,linfarjovTAC14leaderselection,clabuspoo14}, we employ this framework to identify important neurons in the brain of C.\ Elegans.

	\vsp
	
	Three nodes in the network have zero in-degree, i.e., they are not influenced by the rest of the network. Thus, as discussed in Remark~\ref{rem.Sj}, there are three ``leader subsets'', each comprised of one of these nodes. Theorem~\ref{thm.stab} implies that system~\eqref{eq.sys} can only be stable if each of these nodes are leaders.

	\vsp
	
\begin{figure*}[t!] 
\centering
	\begin{tabular}{cccc}
	\!\!\!\!\!
	\!\!\!\!\!
	\subfloat[$J_2$ optimal leaders]{
	\label{fig.celnet2}
	\includegraphics[width=.25\textwidth]{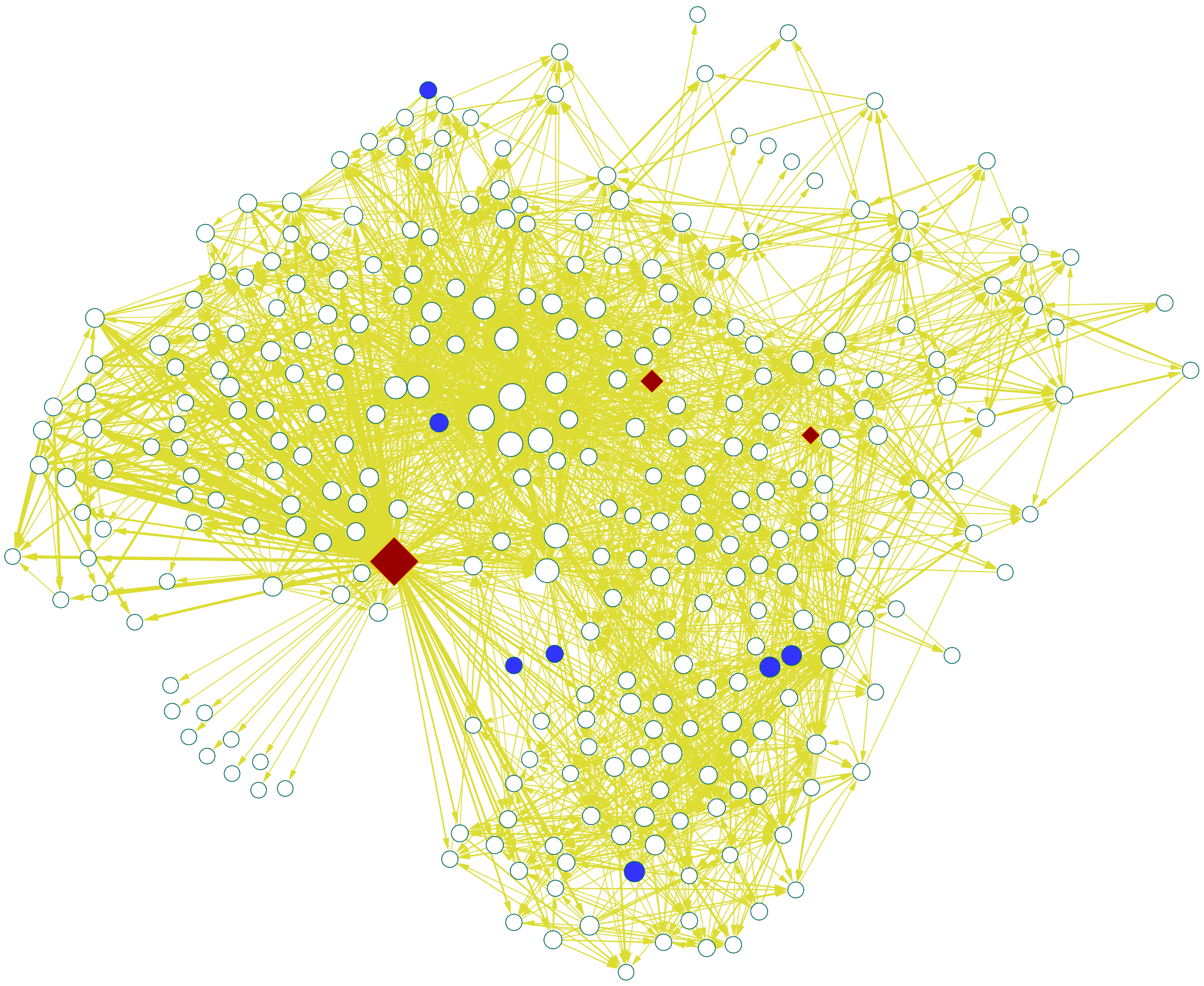}}
	\!\!\!\!\!
	&
	\!\!\!\!\!
	\subfloat[$J_\infty$ optimal leaders]{
	\label{fig.celnetinf}
         \includegraphics[width=.25\textwidth]{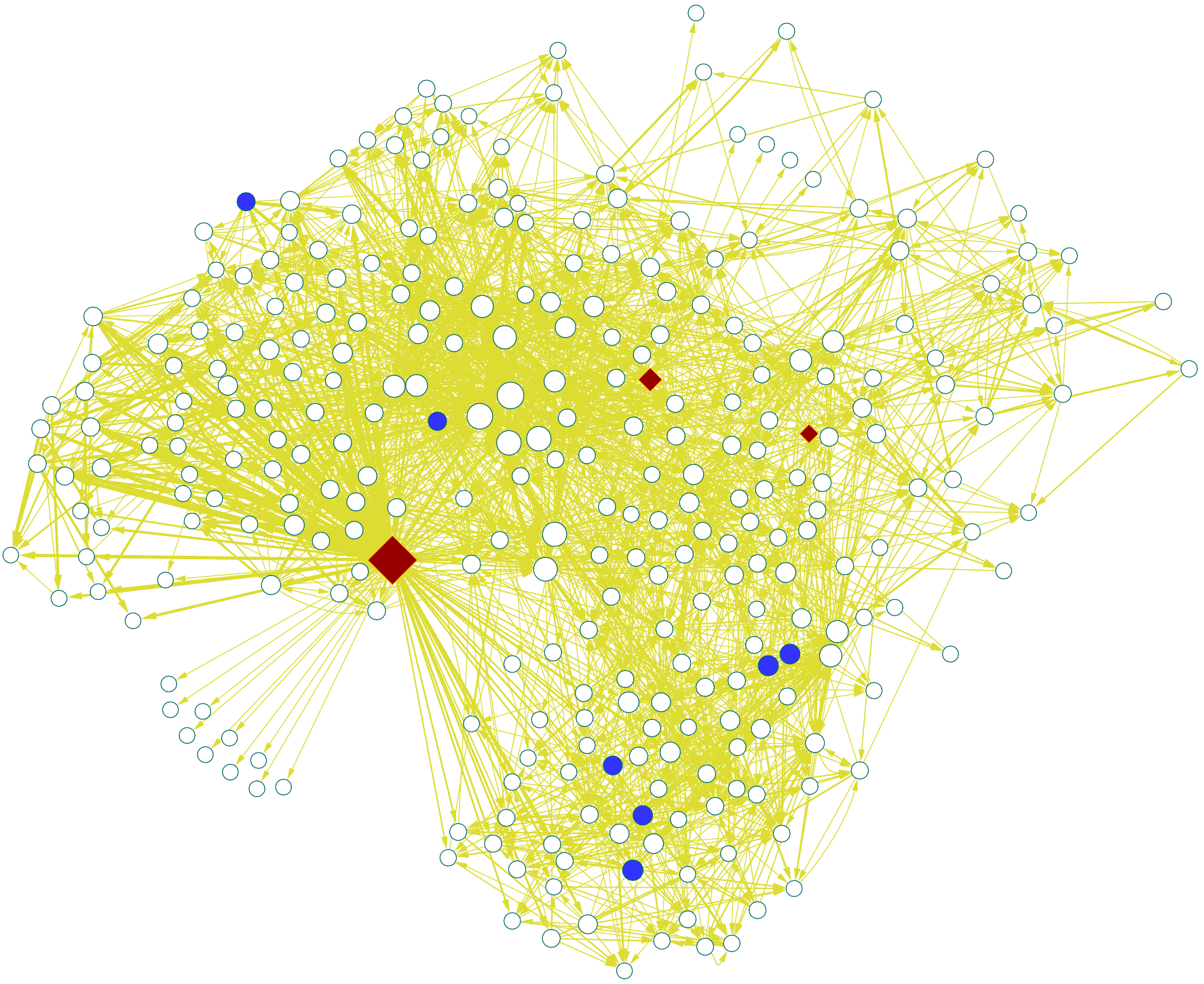}}
	\!\!\!\!\!
         &
	\!\!\!\!\!
         \subfloat[$J_2$ optimal leaders] {
    	\label{fig.celJ2}
	\raisebox{17mm}{
  	\begin{tabular}{rc}
    	\!\!\!\!\!
        \rotatebox{90}{\scriptsize $100*(J(u)/\Jlb_2(N)-1)$}
    	\!\!\!\!\!
    	&
        \!\!\!\!\!
    	\!\!\!\!\!
    	{	 \includegraphics[width=.22\textwidth]{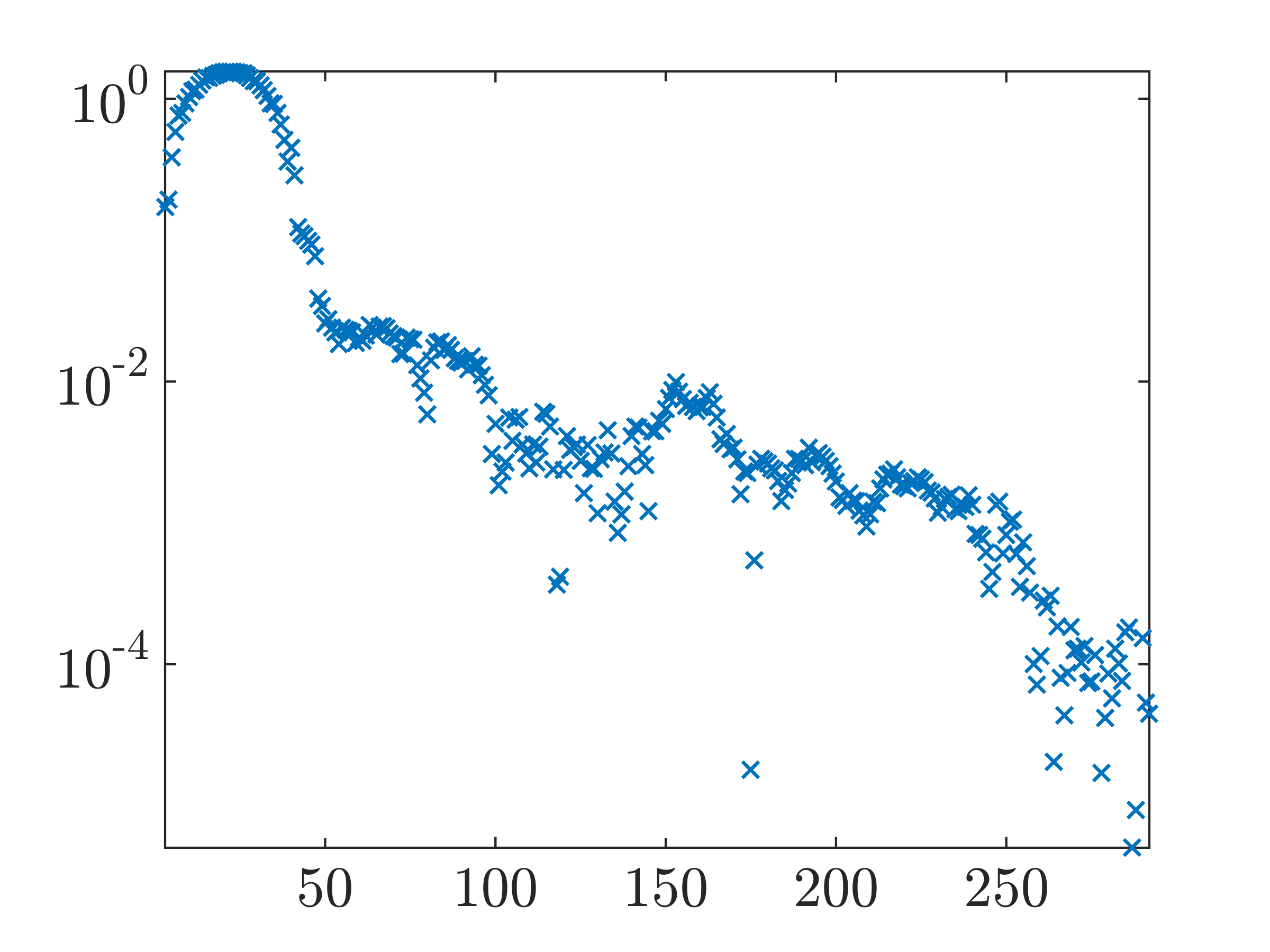}}
    	\!\!\!\!\!    \!\!\!\!\!
    	\\
    	\!\!\!\!\!
        \!\!\!\!\!
    	\!\!\!\!\!
    	&
    	\!\!\!\!\!
    	{\footnotesize Number of leaders $N$}
    	\!\!\!\!\!
    	\end{tabular}}}
	\!\!\!\!\!
	&
	\!\!\!\!\!
	\subfloat[$J_\infty$ optimal leaders] {
    	\label{fig.celJinf}
	\raisebox{17mm}{
  	\begin{tabular}{c}
        \!\!\!\!\!
    	\!\!\!\!\!
    	{	 \includegraphics[width=.22\textwidth]{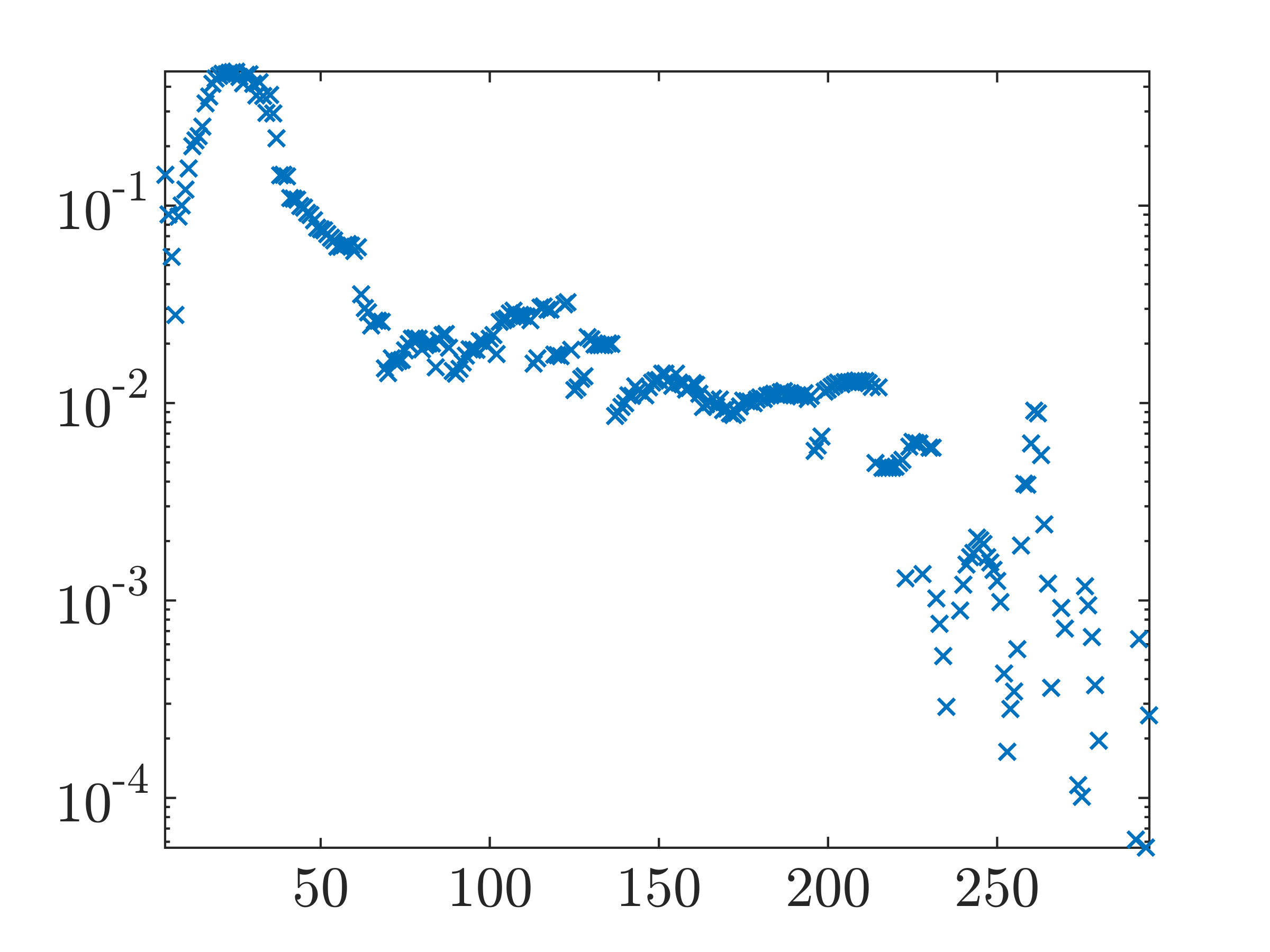}}
    	\!\!\!\!\!    \!\!\!\!\!
    	\\
    	\!\!\!\!\!
    	{\footnotesize 
	Number of leaders $N$}
    	\!\!\!\!\!
    	\end{tabular}}}
	\!\!\!\!\!
	\!\!\!\!\!
        \end{tabular}
        \caption{C.\ Elegans neural network with $N = 10$ (a) $J_2$ and (b) $J_\infty$ leaders  along with the (c) $J_2$ and (d) $J_\infty$ performance of varying numbers of leaders $N$ relative to $\Jlb(N)$. In all cases, leaders are selected via ``rounding''.}
        \label{fig.neth2hi}
\end{figure*}

In Figs.~\ref{fig.celJ2} and~\ref{fig.celJinf}, we show $J_2$ and $J_\infty$ resulting from ``rounding'' the solution to problem~\eqref{pr.rel} to select the additional $1$ to $294$ leaders. Performance is plotted as an increase (in percents) relative to the lower bound $\Jlb(N)$ obtained from~\eqref{pr.rel}. This provides an upper bound on 
suboptimality of the identified set of leaders. While this value does not provide information about how $J(u)$ changes with the number of leaders, Remark~\ref{rem.mono} implies that it monotonically \mbox{decreases with $N$.}

	\vsp
 
For both $J_2$ and $J_\infty$ performance metrics, Figs.~\ref{fig.celJ2} and~\ref{fig.celJinf} illustrates that 
the upper bound is loosest for $25$ leaders ($1.56\%$ and $0.48\%$, respectively). As seen in Fig.~\ref{fig.dirnetperf} from the previous example, whose small size enabled exhaustive search to solve Problem~\ref{pr.comb} exactly, the upper bound on suboptimality is not tight and the exact optimal solution to Problem~\ref{pr.comb} can differ by as much $21.75\%$ from the lower bound. This suggests that ``rounding'' selects very good sets of leaders for this example.

	\vsp

In Figs.~\ref{fig.celnet2} and~\ref{fig.celnetinf}, we show the network with ten identified $J_2$ and $J_\infty$ optimal leaders. The size of the nodes is related to their out-degree and the thickness of the edges is related to the weight. The \tc{cytored}{red $\diamond$} marks nodes that {\em must\/} be leaders and the \tc{cytoblue2}{blue $\circ$} marks the $7$ additional leaders selected by ``rounding''.

	\vspace*{-2ex}
\section{Combination drug therapy}
\label{sec.cdt}

System~\eqref{eq.sys} also arises in the modeling of combination drug therapy~\cite{hernandez2011discrete,jonmatmur13,jonranmur13,jonmatmur14,jon16}, and it provides a model for the evolution of populations of mutants of the HIV virus $x$ in the presence of a combination of drugs $u$. The HIV virus is known to be present in the body in the form of different mutant strands; in~\eqref{eq.sys}, the $i$th component of the state vector $x$ represents the population of the $i$th HIV mutant. The diagonal entries of the matrix $A$ represent the net replication rate of each mutant, and the off diagonal entries of $A$, which are all nonnegative, represent the rate of mutation from one mutant to another. The control input $u_k$ is the dose of drug $k$ and each column $D_k$ of the matrix $D$ in $K(u) = \diag \, (Du)$ specifies at what rate drug $k$ kills each HIV mutant.

	\vspace*{-3ex}
\subsection{Nondifferentiability of $J_\infty$}

\begin{figure}
\[
	\begin{aligned}
	\raisebox{-13mm}{
	\resizebox{.45\hsize}{!}{
	\input{network_ex_drug}}}
	&
	~~~
	&
	\resizebox{.45\hsize}{!}{$
	\underbrace{
         \matbegin \begin{array}{cccc}
                1 & 1 &  & \\
                 1 & 1 &  & \\
                  &  & 1 & 1\\
                 &  & 1 & 1
                \end{array}
                \matend
                }_{A}
                ~
\underbrace{
         \matbegin \begin{array}{cc}
                -1 & 0\\
                 -1 & .1 \\
                  .1 & -1  \\
                  0 & -1 
                \end{array}
                \matend
                }_{K}        
                $}
        \end{aligned}
        \]
\caption{A directed network and corresponding $A$ matrix for a virus with $4$ mutants and $2$ drugs. For this system, $J_\infty$ is nondifferentiable.}
\label{fig.nondif}
\end{figure}
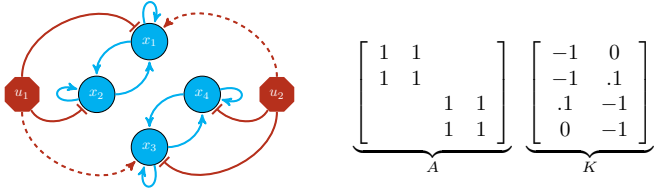

\begin{figure}[t]
\begin{center}
\vspace{-.5cm}
	\begin{tabular}{cc}
	    \!\!\!\!\!
    \!\!\!\!\!
	\subfloat[Descent of objective function.] {
    \label{fig.dirnetperf1}
  \begin{tabular}{rc}
    \!\!\!\!\!
        \rotatebox{90}{\footnotesize \quad\quad\;\; $h(u^k) - h(u^\star)$}
    \!\!\!\!\!
    &
        \!\!\!\!\!
    \!\!\!\!\!
    {	 \includegraphics[width=0.5\columnwidth]{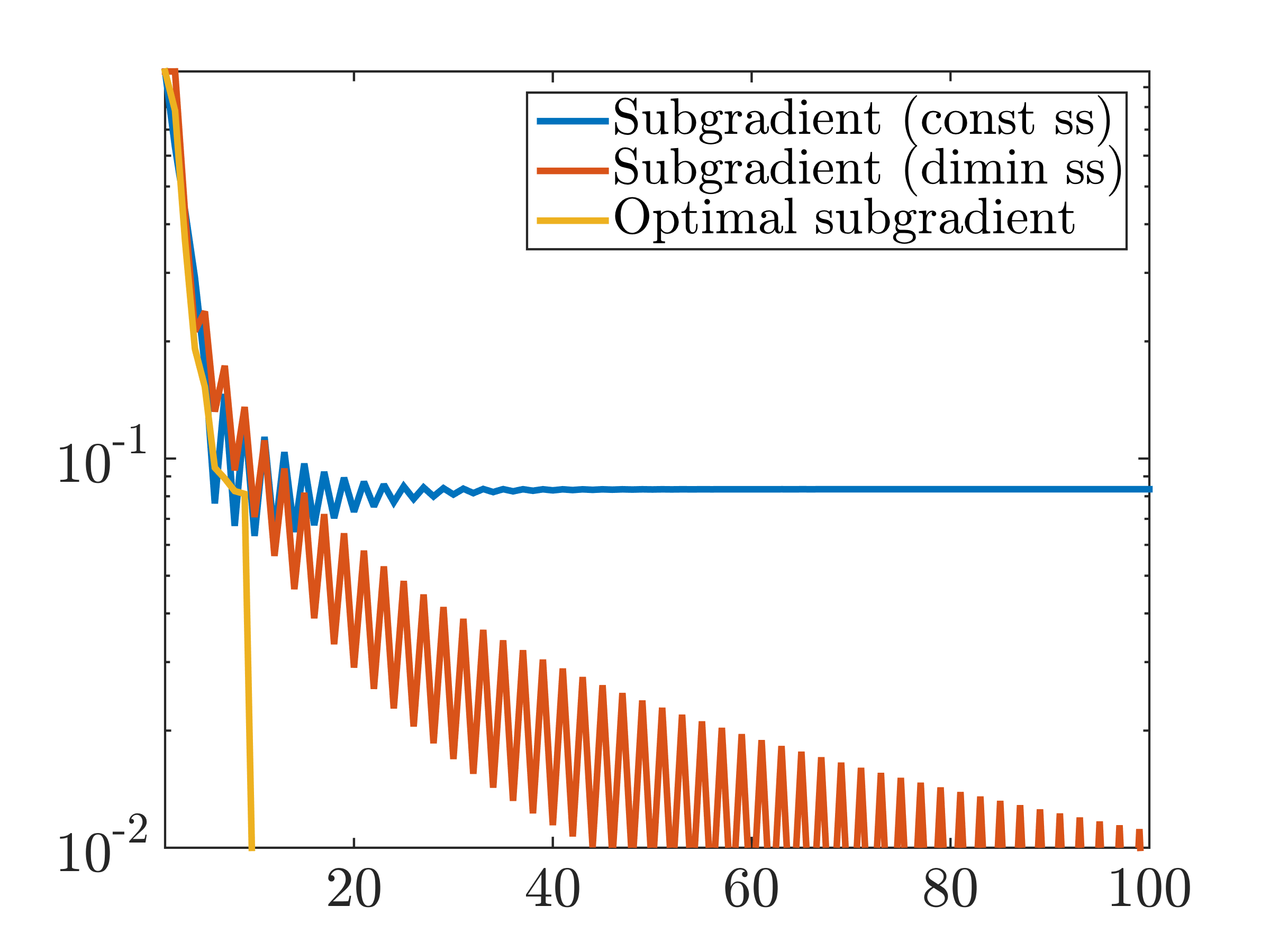}}
    \!\!\!\!\!    \!\!\!\!\!
    \\
    \!\!\!\!\!
        \!\!\!\!\!
    \!\!\!\!\!
    &
    \!\!\!\!\!
    {\footnotesize Iteration $k$}
    \!\!\!\!\!
    \end{tabular}}
        \!\!\!\!\!       \!\!\!\!\!
    &
        \!\!\!\!\!       \!\!\!\!\!
	\subfloat[Iterates  in the ($u_1, u_2)$-plane.] {
    \label{fig.dirnetplane}
  \begin{tabular}{rc}
    \!\!\!\!\!
        \rotatebox{90}{\footnotesize \quad\quad\quad\quad\quad \; $u_2$}
    \!\!\!\!\!
    &
        \!\!\!\!\!
    \!\!\!\!\!
    {	 \includegraphics[width=0.5\columnwidth]{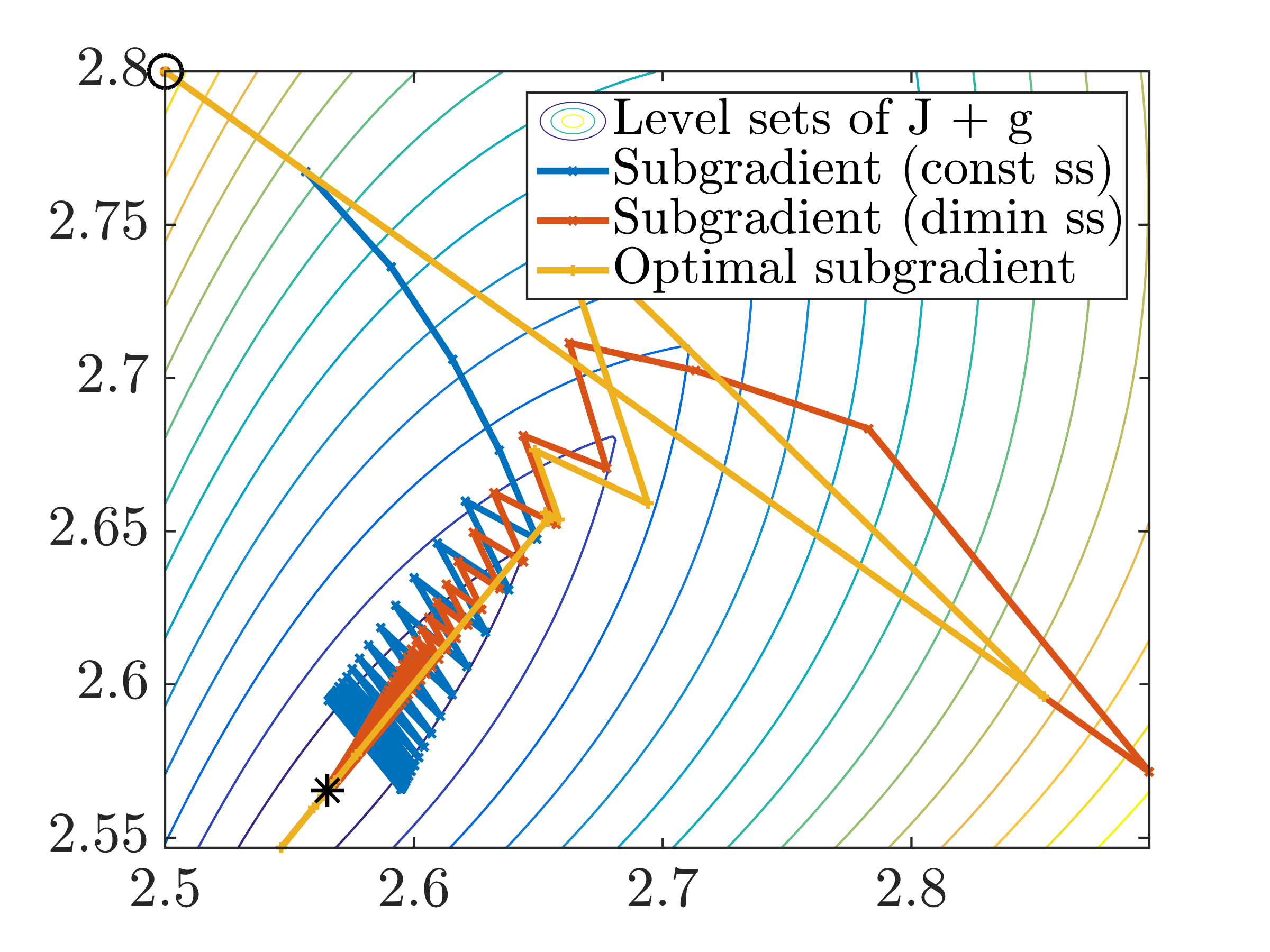}}
    \!\!\!\!\!    \!\!\!\!\!
    \\
    \!\!\!\!\!
        \!\!\!\!\!
    \!\!\!\!\!
    &
    \!\!\!\!\!
    {\footnotesize $u_1$}
    \!\!\!\!\!
    \end{tabular}}
        \!\!\!\!\!
    \!\!\!\!\!
	\end{tabular}
\end{center}
    \caption{Comparison of different algorithms starting from initial condition $[2.5~2.8]^T$. The algorithms are the subgradient method with a constant step-size (\tc{matlabblue}{blue}), the subgradient method with a diminishing step-size (\tc{matlabred}{red}) and our optimal subgradient method where the step-size is chosen via backtracking to ensure descent of the objective function (\tc{matlaborange}{yellow}).
    }
	\label{fig.subg} 
\end{figure}

The mutation patterns of viruses need not be connected. In Fig.~\ref{fig.nondif}, we show a sample mutation network with 2 disconnected components. For this network, the $\cH_\infty$ norm is nondifferentiable when $u_1 = u_2$. Nondifferentiability and the lack of an efficiently computable proximal operator necessitates the use of subgradient methods for solving
\[
	\ds\minimize_u
	~~
	J_\infty(u)
	~+~
	u^Tu.
\]
As shown in Fig.~\ref{fig.subg} with $h (u) \DefinedAs J_\infty(u) + u^Tu$, subgradient methods are not descent methods so small constant or a divergent series of diminishing step-sizes must be employed.

	\vsp
	
We compare the performance of the subgradient method with a constant step-size of $10^{-2}$ (\tc{matlabblue}{blue}) and a diminishing step-size $\tfrac{7 \times 10^{-2}}{k}$ (\tc{matlabred}{red}) with our optimal subgradient method in which the step-size is chosen via backtracking to ensure descent of the objective function (\tc{matlaborange}{yellow}). We show the objective function value with respect to iteration number in Fig.~\ref{fig.dirnetperf1} and the iterates $u^k$ in the $(u_1, u_2)$-plane in Fig.~\ref{fig.dirnetplane}.

	\vsp
	
We run the subgradient methods for $1000$ iterations as there is no principled stopping criterion. Our optimal subgradient method converged with an accuracy of $10^{-4}$ (i.e., there was a $v \in \partial J_\infty(u)$ such that $\norm{v + \nabla g(u)} \leq 10^{-4}$), in $23$ iterations.

	\vspace*{-2ex}
\subsection{A clinically relevant example}

\begin{figure}[t]
\begin{center}
\vspace{-.5cm}
	\begin{tabular}{cc}
	    \!\!\!\!\!
    \!\!\!\!\!
	\subfloat[HIV mutation network]{
    \label{fig.net}
  \begin{tabular}{c}
        \!\!\!\!\!
    \!\!\!\!\!
    {	 \includegraphics[width=0.35\columnwidth]{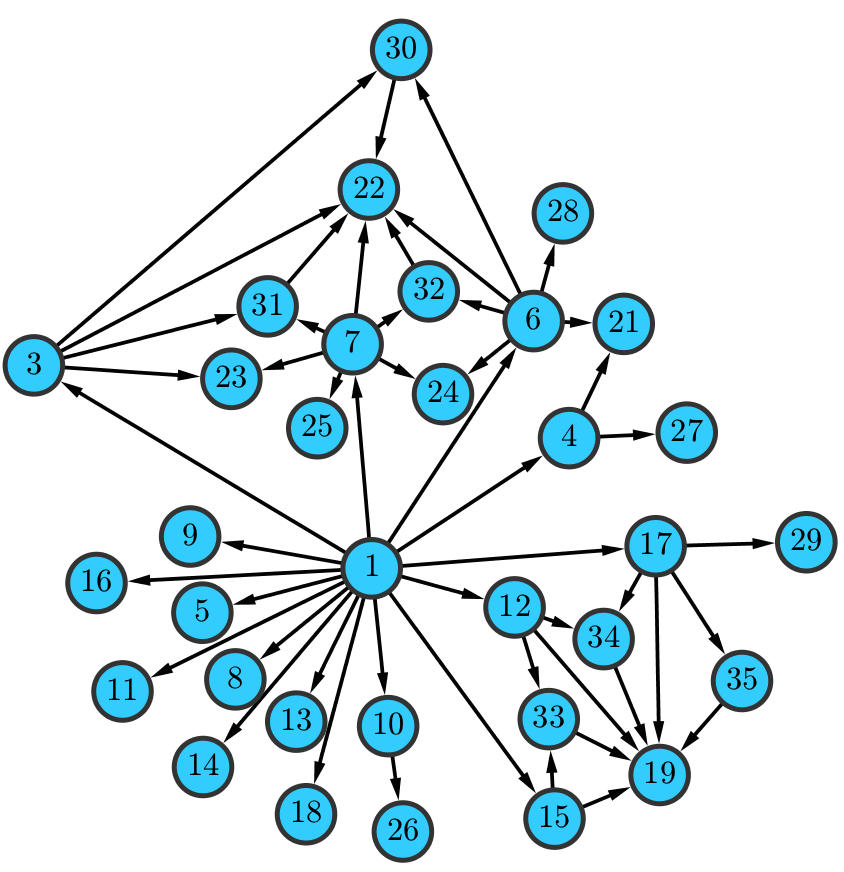}}
    \end{tabular}}
    &
	\subfloat[Sparsity pattern of $A$] {
    \label{fig.sp}
  \begin{tabular}{c}
        \!\!\!\!\!
    \!\!\!\!\!
    {	 \includegraphics[width=0.35\columnwidth]{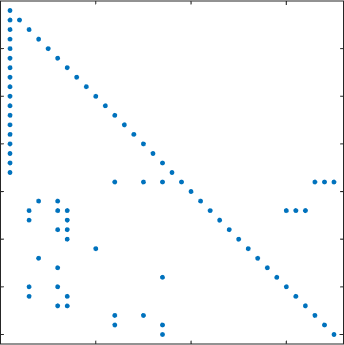}}
    \!\!\!\!\!    \!\!\!\!\!
    \end{tabular}}
        \!\!\!\!\!
    \!\!\!\!\!
	\end{tabular}
\end{center}
    \caption{Mutation pattern of the HIV mutants from~\cite{klehalhor12}.}
	\label{fig.hiv}
\end{figure}

Following~\cite{jonranmur13,jonmatmur13} and using data from~\cite{klehalhor12}, we study a system with $35$ mutants and $5$ drugs. The sparsity pattern of the matrix $A$, shown in Fig.~\ref{fig.hiv}, corresponds to the mutation pattern and replication rates  of $33$ mutants and $K(u)$ specifies the effect of drug therapy. Two mutants are not shown in Fig.~\ref{fig.net} as they have no mutation pathways to or from other mutants.

	\vsp

Several clinically relevant properties, such as maximum dose or budget constraints, may be directly enforced in our formulation. Other combinatorial conditions can be promoted via convex penalties, such as drug $j$ requiring drug $i$ via $u_i \geq u_j$ or mutual exclusivity of drugs $i$ and $j$ via $u_i + u_j \leq 1$. We design optimal drug doses \mbox{using two convex regularizers $g$.}

	\vsp

\subsubsection{Budget constraint}
We impose a unit budget constraint on the drug doses and solve the $J_2$ and $J_\infty$ problems using proximal gradient methods~\cite{becteb09,ber99}. These can be cast in the form~\eqref{pr.gen}, where $g$ is the indicator function associated with the probability simplex $\cal P$. Table~\ref{tab.dos} contains the optimal doses and illustrates the tradeoff between $\Ht$ and $\Hinf$ performance.

\begin{table}[htb]
\centering
\begin{tabular}{|r|c|c|}
    \hline
    \bf {Antibody} & $\tc{matlabred}{\bf u^\star_{2}}$ & $\tc{matlabblue}{u^\star_{\infty}}$ \\ \hline
    3BC176    &    \tc{matlabred}{$0.5952$}&    \tc{matlabblue}{$0.9875$} \\ \hline
    PG16      &    \tc{matlabred}{$0$}&         \tc{matlabblue}{0}          \\ \hline
    45-46G54W &    \tc{matlabred}{$0.2484$}&    \tc{matlabblue}{$0.0125$} \\ \hline
    PGT128    &    \tc{matlabred}{$0.1564$}&    \tc{matlabblue}{$0$}     \\ \hline
    10-1074   &    \tc{matlabred}{$0$}&    \tc{matlabblue}{$0$}          \\ \hline
\end{tabular}
~
\begin{tabular}{|c|c|}
    \hline
    $J_2(\tc{matlabred}{u^\star_{2}})$ & $0.6017$ \\ \hline
    $J_2(\tc{matlabblue}{u^\star_{\infty}})$ & $1.1947$ \\ \hline \hline
    $J_\infty(\tc{matlabred}{u^\star_{2}})$ & $0.1857$ \\ \hline
    $J_\infty(\tc{matlabblue}{u^\star_{\infty}})$ & $0.1084$ \\ \hline
\end{tabular}
\caption{Optimal budgeted doses and $\cH_2$/$\cH_\infty$ performance.}
\label{tab.dos}
\end{table}

	\vsp
	
\subsubsection{Sparsity-promoting framework}
\begin{figure}
      \centering
      \begin{tabular}{rc}
    \!\!\!\!\!
        \rotatebox{90}{\footnotesize \quad \quad Percent}
        \rotatebox{90}{\footnotesize degradation of $J$}
    \!\!\!\!\!    \!\!\!
    &
    \!\!\!\!\!    
    {\includegraphics[width=0.9\columnwidth]{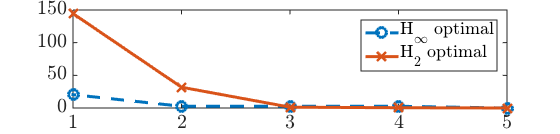}}
    \!\!\!\!\!    \!\!\!\!\!
    \\
    \!\!\!\!\!
        \!\!\!\!\!
    \!\!\!\!\!
    &
    \!\!\!\!\!
    {\footnotesize Number of antibodies}
    \!\!\!\!\!
    \end{tabular}
      \caption{Performance degradation (in percents) relative to the optimal \tc{matlabred}{$\cH_2$} and \tc{matlabblue}{$\cH_\infty$} strategies that use all $5$ drugs.}
      \label{fig.spprom}
\end{figure}

Although the above budget constraint is naturally sparsity-promoting, in Algorithm~\ref{alg:HIV} we augment a quadratically regularized optimal control problem with a reweighted $\ell_1$ norm~\cite{canwakboy08} to select a homotopy path of successively sparser sets of drugs. We then perform a `polishing' step to design the optimal doses of the selected set of drugs. We use $50$ logarithmically spaced increments of the regularization parameter $\gamma$ between $0.01$ and $10$ to identify the drugs and then replace the weighted $\ell_1$ penalty with a constraint to prescribe the selected drugs. In Fig.~\ref{fig.spprom}, we show performance degradation (in percents) relative to the optimal dose that uses all $5$ drugs with $B = C = I$ and $R = I$.

\begin{algorithm}
 Set $\gamma>0$, $R \succ0$, $w = \mathds{1}$, $\varepsilon >0$

 \While{$\card(u_\gamma)>N$}{
 $
	u_\gamma \; = \; \argmin_u ~ J(u) \,+\, u^T R \, u \,+\, \gamma \, \sum_i w_i \, |u_i |
 $
 
 increase $\gamma$,~ set $w_i \, = \, 1/(u_i \,+\, \varepsilon)$
 }
 $
 \ba{rrl}
	u^\star_N
	\; = \;
	&
	\!\!
	\argmin_u & J(u) \;+\; u^T R \, u
	\\
	&
	\!\!
	\subject 
	&
	 \sparse(u) \, \subseteq \, \sparse(u_\gamma)
	\ea
 $
 \caption{Sparsity-promoting algorithm for $N$ drugs}
 \label{alg:HIV}
\end{algorithm}

	\vspace*{-1ex}
\section{Concluding remarks}
\label{sec:conc}

We introduce a unifying framework for the $\cH_2$ and $\cH_\infty$ synthesis of positive systems and use it to address the problems of leader selection in directed consensus networks and the design of combination drug therapy for HIV treatment. We identify classes of networks for which the $\cH_\infty$ norm is a differentiable function of the control input and develop efficient customized algorithms that perform well even in the absence of differentiability. Our ongoing work focuses on the design of time-varying strategies within an MPC framework.

	\vspace*{-1ex}
\section*{Acknowledgments}
\label{sec.ack}

We would like to thank Anders Rantzer and Katie Fitch for useful discussion on positive systems and leader selection, respectively, and to Vanessa Jonsson for providing us with the HIV model and for her insights on combination drug therapy.

\end{document}

%% file: commands.tex
\newtheorem{theorem}{Theorem}
\newtheorem{corollary}[theorem]{Corollary}
\newtheorem{lemma}[theorem]{Lemma}
\newtheorem{proposition}[theorem]{Proposition}
\newtheorem{definition}{Definition}
\newtheorem{remark}{Remark}
\newtheorem{problem}{Problem}

\newcommand{\ds}{\displaystyle}

\newcommand{\enma}[1]   {\ensuremath{#1}}

\newcommand{\non}{\nonumber}

\newcommand{\beq}{\begin{equation}}
\newcommand{\eeq}{\end{equation}}
\newcommand{\bseq}{\begin{subequations}}
\newcommand{\eseq}{\end{subequations}}
\newcommand{\beqn}{\begin{eqnarray}}
\newcommand{\eeqn}{\end{eqnarray}}
\newcommand{\ba}{\begin{array}}
\newcommand{\ea}{\end{array}}
\newcommand{\bct}{\begin{center}}
\newcommand{\ect}{\end{center}}
\newcommand{\btmz}{\begin{itemize}}
\newcommand{\etmz}{\end{itemize}}
\newcommand{\benum}{\begin{enumerate}}
\newcommand{\eenum}{\end{enumerate}}

\newcommand{\mc}{\mathcal}

\newcommand{\R}{{\mathbb R}}

\newcommand{\cH}{\enma{\mathcal H}}
\newcommand{\cL}{\enma{\mathcal L}}

\newcommand{\norm}[1]{\| #1 \|}                 

\newcommand{\Hinf}{\mathcal{H}_{\infty} }
\newcommand{\Ht}{ \mathcal{H}_{2}}

\newcommand{\diag}      {\enma{\mathrm{diag}}}

\newcommand{\trace}     {\enma{\mathrm{trace}}}

\newcommand{\inner}[2]{\left\langle #1,#2 \right\rangle}

\newcommand{\matbegin}{
        \left[
}
\newcommand{\matend}{
        \right]
}

\newcommand{\be}{\begin{equation}}
\newcommand{\ee}{\end{equation}}

\newcommand{\cplxs}{ C\kern -.35em \rule{0.03 em}{.7 ex}~   }

\def\complex{\hbox{C\kern -.45em \rule{0.03 em}{1.5 ex}}~}

\newcommand{\bi}{\begin{itemize}}
\newcommand{\ei}{\end{itemize}}

\newtheorem{assumption}{Assumption}



%% file: smallnetwork.tex
\begin{tikzpicture}[>=stealth',shorten >=1pt, node distance=1.5cm, on grid, initial/.style={}]

  \node[state]          (1)                        {$1$};
  \node[state]          (2) [below right =of 1]    {$2$};
  \node[state]          (3) [above right =of 2]    {$3$};
  \node[state]          (4) [below left =of 2]    {$4$};

\tikzset{every node/.style={fill=white}}

\tikzset{mystyle/.style={->,double=black}}
\path (2)     edge [mystyle]    (3)
      (2)     edge [mystyle]     (4)
      (1)     edge [mystyle]     (2)
      (2)     edge [mystyle]     (1);
      
\end{tikzpicture} 

%% file: dirnet_bal_med.tex
\begin{tikzpicture}[>=stealth',shorten >=1pt, node distance=1.1cm, on grid, initial/.style={}]

  \node[state,scale=0.7,fill=matlabred,text=white]          (8)     {$8$};
  \node[state,scale=0.7,fill=matlabblue,text=white]          (7) [below =of 8]    {$7$};
  \node[state,scale=0.7]          (2) [right =of 7]    {$2$};
  \node[state,scale=0.7]          (1) [right =of 8]                       {$1$};
  \node[state,scale=0.7]          (3) [below =of 7]    {$3$};
  \node[state,scale=0.7]          (6) [above =of 8]    {$6$};
  \node[state,scale=0.7,fill=matlaborange]          (4) [left =of 8]    {$4$};
  \node[state,scale=0.7]          (5) [left =of 7]    {$5$};

\tikzset{every node/.style={fill=white}}

\tikzset{mystyle/.style={->,double=black}}
\path (7)     edge [mystyle]   (5)
(4)     edge [mystyle]    (5)
      (8)     edge [mystyle]    (6)
      (2)     edge [mystyle] node[below=0.1cm,scale=1] {$2$}    (7)
      (1)     edge [mystyle]    (2);

%


\tikzset{mystyle/.style={->,relative=false,in=-90,out=0,double=black}}
\path (3)     edge [mystyle]    (2);

\tikzset{mystyle/.style={->,relative=false,in=-180,out=-90,double=black}}
\path (5)     edge [mystyle]    (3);

\tikzset{mystyle/.style={->,relative=false,in=-135,out=45,double=black}}
\path (5)     edge [mystyle]    (8);

\tikzset{mystyle/.style={->,relative=false,in=90,out=0,double=black}}
\path (6)     edge [mystyle]   (1);

\tikzset{mystyle/.style={->,relative=false,in=-70,out=70,double=black}}
\path (7)     edge [mystyle] node[right=0.1cm,scale=1] {$2$}  (8);

\tikzset{mystyle/.style={->,relative=false,in=110,out=-110,double=black}}
\path (8)     edge [mystyle]   (7);


\tikzset{mystyle/.style={->,relative=false,in=0,out=180,double=black}}
\path (8)     edge [mystyle]   (4);

\end{tikzpicture} 

%% file: network_ex_drug.tex
\begin{tikzpicture}[>=stealth',shorten >=1pt, node distance=1.8cm, on grid, initial/.style={}]

  \node[state,fill=ProcessBlue,text=white]          (1)                        {$x_1$};
  \node[state,fill=ProcessBlue,text=white]          (2) [below left =of 1]    {$x_2$};
  \node[state,fill=ProcessBlue,text=white]          (4) [below right =of 1]    {$x_4$};
  \node[state,fill=ProcessBlue,text=white]          (3) [below left =of 4]    {$x_3$};
  \node[regular polygon, regular polygon sides=8,fill=BrickRed,text=white]          (5) [right =of 4]    {$u_2$};
  \node[regular polygon, regular polygon sides=8,fill=BrickRed,text=white]          (6) [left =of 2]    {$u_1$};

\tikzset{every node/.style={fill=white}}


%

\tikzset{mystyle/.style={->,relative=false,in=30,out=150,double=ProcessBlue,ProcessBlue,loop above,min distance = 20}}
\path (1)     edge [mystyle]   (1);

\tikzset{mystyle/.style={->,relative=false,in=70,out=110,double=ProcessBlue,ProcessBlue,loop left,min distance = 20}}
\path (2)     edge [mystyle]   (2);

\tikzset{mystyle/.style={->,relative=false,in=70,out=110,double=ProcessBlue,ProcessBlue,loop below,min distance = 20}}
\path (3)     edge [mystyle]   (3);

\tikzset{mystyle/.style={->,relative=false,in=70,out=110,double=ProcessBlue,ProcessBlue,loop right,min distance = 20}}
\path (4)     edge [mystyle]   (4);

%

%

\tikzset{mystyle/.style={->,relative=false,in=-90,out=0,double=ProcessBlue,ProcessBlue}}
\path (3)     edge [mystyle]   (4);

\tikzset{mystyle/.style={->,relative=false,in=90,out=180,double=ProcessBlue,ProcessBlue}}
\path (4)     edge [mystyle]   (3);

\tikzset{mystyle/.style={->,relative=false,in=90,out=180,double=ProcessBlue,ProcessBlue}}
\path (1)     edge [mystyle]   (2);

\tikzset{mystyle/.style={->,relative=false,in=-90,out=0,double=ProcessBlue,ProcessBlue}}
\path (2)     edge [mystyle]   (1);

\tikzset{mystyle/.style={->,relative=false,in=45,out=90,double=BrickRed,BrickRed,dashed}}
\path (5)     edge [mystyle]   (1);

\tikzset{mystyle/.style={-|,relative=false,in=-45,out=-90,double=BrickRed,BrickRed}}
\path (5)     edge [mystyle]   (3);

\tikzset{mystyle/.style={-|,relative=false,in=-45,out=-120,double=BrickRed,BrickRed}}
\path (5)     edge [mystyle]   (4);

\tikzset{mystyle/.style={-|,relative=false,in=135,out=90,double=BrickRed,BrickRed}}
\path (6)     edge [mystyle]   (1);

\tikzset{mystyle/.style={-|,relative=false,in=-135,out=-60,double=BrickRed,BrickRed}}
\path (6)     edge [mystyle]   (2);

\tikzset{mystyle/.style={->,relative=false,in=-135,out=-90,double=BrickRed,BrickRed,dashed}}
\path (6)     edge [mystyle]   (3);

\end{tikzpicture} 

%% file: decentr-revised-arXiv.bbl
\begin{thebibliography}{10}

\bibitem{bamvou05}
B.~Bamieh and P.~G. Voulgaris, ``A convex characterization of distributed
  control problems in spatially invariant systems with communication
  constraints,'' {\em Syst. and Control Lett.}, vol.~54, pp.~575--583, 2005.

\bibitem{rotlalTAC06}
M.~Rotkowitz and S.~Lall, ``A characterization of convex problems in
  decentralized control,'' {\em IEEE Trans. Automat. Control}, vol.~51, no.~2,
  pp.~274--286, 2006.

\bibitem{wang2016system}
Y.-S. Wang, N.~Matni, and J.~C. Doyle, ``A system level approach to controller
  synthesis,'' {\em IEEE Trans. Automat. Control}, 2016.
\newblock submitted; also arXiv:1610.04815.

\bibitem{farrin11}
L.~Farina and S.~Rinaldi, {\em Positive Linear Systems: Theory and
  Applications}.
\newblock John Wiley \& Sons, 2011.

\bibitem{hernandez2011discrete}
E.~Hernandez-Vargas, P.~Colaneri, R.~Middleton, and F.~Blanchini,
  ``Discrete-time control for switched positive systems with application to
  mitigating viral escape,'' {\em Int. J. Robust Nonlin.}, vol.~21, no.~10,
  pp.~1093--1111, 2011.

\bibitem{herestmid13}
E.~A. Hernandez-Vargas and R.~H. Middleton, ``Modeling the three stages in
  {HIV} infection,'' {\em J. Theor. Biol.}, vol.~320, pp.~33--40, 2013.

\bibitem{jonmatmur13}
V.~Jonsson, N.~Matni, and R.~M. Murray, ``Reverse engineering combination
  therapies for evolutionary dynamics of disease: An {${\cal H}_\infty$}
  approach,'' in {\em Proceedings of the 52nd IEEE Conference on Decision and
  Control}, pp.~2060--2065, 2013.

\bibitem{jonranmur13}
V.~Jonsson, A.~Rantzer, and R.~M. Murray, ``A scalable formulation for
  engineering combination therapies for evolutionary dynamics of disease,''
  {\em Proceedings of the 2014 American Control Conference}, pp.~2771--2778,
  2014.

\bibitem{jonmatmur14}
V.~Jonsson, N.~Matni, and R.~M. Murray, ``Synthesizing combination therapies
  for evolutionary dynamics of disease for nonlinear pharmacodynamics,'' in
  {\em Proceedings of the 53rd Conference on Decision and Control},
  pp.~2352--2358, 2014.

\bibitem{jon16}
V.~D. Jonsson, {\em Robust control of evolutionary dynamics}.
\newblock PhD thesis, California Institute of Technology, 2016.

\bibitem{patbam10}
S.~Patterson and B.~Bamieh, ``Leader selection for optimal network coherence,''
  in {\em Proceedings of the 49th IEEE Conference on Decision and Control},
  pp.~2692--2697, 2010.

\bibitem{linfarjovTAC14leaderselection}
F.~Lin, M.~Fardad, and M.~R. Jovanovi\'c, ``Algorithms for leader selection in
  stochastically forced consensus networks,'' {\em IEEE Trans. Automat.
  Control}, vol.~59, pp.~1789--1802, July 2014.

\bibitem{fitch2013information}
K.~Fitch and N.~E. Leonard, ``Information centrality and optimal leader
  selection in noisy networks,'' in {\em Proceedings of the 52nd IEEE
  Conference on Decision and Control}, pp.~7510--7515, 2013.

\bibitem{fitch2014joint}
K.~Fitch and N.~Leonard, ``Joint centrality distinguishes optimal leaders in
  noisy networks,'' {\em IEEE Trans. Control Netw. Syst.}, vol.~3, no.~4,
  pp.~366--378, 2016.

\bibitem{clabuspoo14}
A.~Clark, L.~Bushnell, and R.~Poovendran, ``A supermodular optimization
  framework for leader selection under link noise in linear multi-agent
  systems,'' {\em IEEE Trans. Automat. Control}, vol.~59, no.~2, pp.~283--296,
  2014.

\bibitem{clark2016submodularity}
A.~Clark, B.~Alomair, L.~Bushnell, and R.~Poovendran, {\em Submodularity in
  dynamics and control of networked systems}.
\newblock Springer, 2016.

\bibitem{tanlan11}
T.~Tanaka and C.~Langbort, ``The bounded real lemma for internally positive
  systems and $\mathcal{H}_\infty$ structured static state feedback,'' {\em
  IEEE Trans. Automat. Control}, vol.~56, no.~9, pp.~2218--2223, 2011.

\bibitem{ran15}
A.~Rantzer, ``Scalable control of positive systems,'' {\em Eur. J. Control},
  vol.~24, pp.~72--80, 2015.

\bibitem{bri13}
C.~Briat, ``Robust stability and stabilization of uncertain linear positive
  systems via integral linear constraints: $\mathcal{L}_1$ gain and
  $\mathcal{L}_\infty$ gain characterization,'' {\em Int. J. Robust Nonlin.},
  vol.~23, no.~17, pp.~1932--1954, 2013.

\bibitem{ebipeaarz11}
Y.~Ebihara, D.~Peaucelle, and D.~Arzelier, ``$\mathcal{L}_1$ gain analysis of
  linear positive systems and its application,'' in {\em Proceedings of 50th
  IEEE Conference on Decision and Control}, pp.~4029--4034, 2011.

\bibitem{ColSmi:2014:IFA_4769}
M.~Colombino and R.~Smith, ``Convex characterization of robust stability
  analysis and control synthesis for positive linear systems,'' in {\em
  Proceedings of the 53rd IEEE Conference on Decision and Control},
  pp.~4379--4384, 2014.

\bibitem{ColSmi:2016:IFA_5242}
M.~Colombino and R.~S. Smith, ``A convex characterization of robust stability
  for positive and positively dominated linear systems,'' {\em IEEE Trans.
  Automat. Control.}, vol.~61, no.~7, pp.~1965--1971, 2016.

\bibitem{ran16}
A.~Rantzer, ``On the {K}alman-{Y}akubovich-{P}opov lemma for positive
  systems,'' {\em IEEE Trans. Automat. Control}, vol.~61, no.~5,
  pp.~1346--1349, 2016.

\bibitem{colaneri2014convexity}
P.~Colaneri, R.~H. Middleton, Z.~Chen, D.~Caporale, and F.~Blanchini,
  ``Convexity of the cost functional in an optimal control problem for a class
  of positive switched systems,'' {\em Automatica}, vol.~4, no.~50,
  pp.~1227--1234, 2014.

\bibitem{ranber14}
A.~Rantzer and B.~Bernhardsson, ``Control of convex monotone systems,'' in {\em
  Proceedings of the 53rd IEEE Conference on Decision and Control},
  pp.~2378--2383, 2014.

\bibitem{herestcol14}
E.~A. Hernandez-Vargas, P.~Colaneri, and R.~H. Middleton, ``Switching
  strategies to mitigate hiv mutation,'' {\em IEEE Trans. Control Syst.
  Technol.}, vol.~22, no.~4, pp.~1623--1628, 2014.

\bibitem{blacolval15}
F.~Blanchini, P.~Colaneri, M.~E. Valcher, {\em et~al.}, ``Switched positive
  linear systems,'' {\em Foundations and {T}rends{\textregistered} in {S}ystems
  and {C}ontrol}, vol.~2, no.~2, pp.~101--273, 2015.

\bibitem{colmidbla16}
P.~Colaneri, R.~H. Middleton, and F.~Blanchini, ``Optimal control of a class of
  positive {M}arkovian bilinear systems,'' {\em Nonlinear {A}nalysis: {H}ybrid
  {S}ystems}, vol.~21, pp.~155--170, 2016.

\bibitem{berple94}
A.~Berman and R.~J. Plemmons, {\em Nonnegative matrices in the mathematical
  sciences}, vol.~9.
\newblock {SIAM}, 1994.

\bibitem{horn1985matrix}
R.~A. Horn and C.~R. Johnson, {\em Matrix Analysis}.
\newblock Cambridge {U}niversity {P}ress, 1985.

\bibitem{dulpag00}
G.~E. Dullerud and F.~Paganini, {\em A course in robust control theory}.
\newblock Springer, 2000.

\bibitem{linfarjovTAC13admm}
F.~Lin, M.~Fardad, and M.~R. Jovanovi\'c, ``Design of optimal sparse feedback
  gains via the alternating direction method of multipliers,'' {\em IEEE Trans.
  Automat. Control}, vol.~58, pp.~2426--2431, September 2013.

\bibitem{jovdhiEJC16}
M.~R. Jovanovi\'c and N.~K. Dhingra, ``Controller architectures: tradeoffs
  between performance and structure,'' {\em Eur. J. Control}, vol.~30,
  pp.~76--91, July 2016.

\bibitem{cohen1981convexity}
J.~E. Cohen, ``Convexity of the dominant eigenvalue of an essentially
  nonnegative matrix,'' {\em Proceedings of the American Mathematical Society},
  vol.~81, no.~4, pp.~657--658, 1981.

\bibitem{boyvan04}
S.~Boyd and L.~Vandenberghe, {\em Convex {O}ptimization}.
\newblock Cambridge {U}niversity {P}ress, 2004.

\bibitem{shor2012minimization}
N.~Z. Shor, {\em Minimization {M}ethods for {N}on-{D}ifferentiable
  {F}unctions}, vol.~3.
\newblock Springer Science \& Business Media, 2012.

\bibitem{bullo2009distributed}
F.~Bullo, J.~Cort{\'e}s, and S.~Martinez, {\em Distributed Control of Robotic
  Networks: A Mathematical Approach to Motion Coordination Algorithms}.
\newblock Princeton University Press, 2009.

\bibitem{becteb09}
A.~Beck and M.~Teboulle, ``A fast iterative shrinkage-thresholding algorithm
  for linear inverse problems,'' {\em SIAM {J.} {I}maging {S}ci.}, vol.~2,
  no.~1, pp.~183--202, 2009.

\bibitem{boyparchupeleck11}
S.~Boyd, N.~Parikh, E.~Chu, B.~Peleato, and J.~Eckstein, ``Distributed
  optimization and statistical learning via the alternating direction method of
  multipliers,'' {\em Found. Trends Mach. Learning}, vol.~3, no.~1, pp.~1--124,
  2011.

\bibitem{ber99}
D.~P. Bertsekas, {\em Nonlinear programming}.
\newblock Athena {S}cientific, 1999.

\bibitem{parboy13}
N.~Parikh and S.~Boyd, ``Proximal algorithms,'' {\em Foundations and Trends in
  optimization}, vol.~1, no.~3, pp.~123--231, 2013.

\bibitem{dhikhojovTAC16}
N.~K. Dhingra, S.~Z. Khong, and M.~R. Jovanovi\'c, ``The proximal augmented
  {L}agrangian method for nonsmooth composite optimization,'' {\em IEEE Trans.
  Automat. Control}, 2016.
\newblock submitted; also arXiv:1610.04514.

\bibitem{xiaboykim07}
L.~Xiao, S.~Boyd, and S.-J. Kim, ``Distributed average consensus with
  least-mean-square deviation,'' {\em J. Parallel Distrib. Comput.}, vol.~67,
  no.~1, pp.~33--46, 2007.

\bibitem{cvetkovic1980spectra}
D.~M. Cvetkovi{\'c}, M.~Doob, and H.~Sachs, {\em Spectra of {G}raphs: {T}heory
  and {A}pplication}.
\newblock Academic Press, 1980.

\bibitem{bamjovmitpat12}
B.~Bamieh, M.~R. Jovanovi\'c, P.~Mitra, and S.~Patterson, ``Coherence in
  large-scale networks: dimension dependent limitations of local feedback,''
  {\em IEEE Trans. Automat. Control}, vol.~57, pp.~2235--2249, September 2012.

\bibitem{lawwoo66}
E.~L. Lawler and D.~E. Wood, ``Branch-and-bound methods: a survey,'' {\em
  Operations Research}, vol.~14, no.~4, pp.~699--719, 1966.

\bibitem{wancan15}
W.~Wang and C.~Lu, ``Projection onto the capped simplex,'' {\em arXiv preprint
  arXiv:1503.01002}, 2015.

\bibitem{dhijovACC15}
N.~K. Dhingra and M.~R. Jovanovi\'c, ``Convex synthesis of symmetric
  modifications to linear systems,'' in {\em Proceedings of the 2015 American
  Control Conference}, pp.~3583--3588, 2015.

\bibitem{dhiwujovNOL16}
N.~K. Dhingra, X.~Wu, and M.~R. Jovanovi\'c, ``Sparsity-promoting optimal
  control of systems with invariances and symmetries,'' in {\em Proceedings of
  the 10th IFAC Symposium on Nonlinear Control Systems}, pp.~648--653, 2016.

\bibitem{zelschall13}
D.~Zelazo, S.~Schuler, and F.~Allg{\"o}wer, ``Performance and design of cycles
  in consensus networks,'' {\em Syst. Control Lett.}, vol.~62, no.~1,
  pp.~85--96, 2013.

\bibitem{watstr98}
D.~J. Watts and S.~H. Strogatz, ``Collective dynamics of `small-world'
  networks,'' {\em Nature}, vol.~393, no.~6684, p.~440, 1998.

\bibitem{whisouthobre86}
J.~G. White, E.~Southgate, J.~N. Thomson, and S.~Brenner, ``The structure of
  the nervous system of the nematode {C}aenorhabditis {E}legans: the mind of a
  worm,'' {\em Phil. Trans. R. Soc. Lond}, vol.~314, pp.~1--340, 1986.

\bibitem{klehalhor12}
F.~Klein, A.~Halper-Stromberg, J.~A. Horwitz, H.~Gruell, J.~F. Scheid,
  S.~Bournazos, H.~Mouquet, L.~A. Spatz, R.~Diskin, A.~Abadir, {\em et~al.},
  ``{HIV} therapy by a combination of broadly neutralizing antibodies in
  humanized mice,'' {\em Nature}, vol.~492, no.~7427, pp.~118--122, 2012.

\bibitem{canwakboy08}
E.~J. Cand\`{e}s, M.~B. Wakin, and S.~P. Boyd, ``Enhancing sparsity by
  reweighted $\ell_1$ minimization,'' {\em J. Fourier Anal. Appl.}, vol.~14,
  pp.~877--905, 2008.

\end{thebibliography}


\vspace*{-1ex}
